\documentclass[11pt]{article}
\usepackage{amssymb,amsmath,latexsym}

\usepackage{color}

\usepackage{epsfig}

\newtheorem{thm}{Theorem}[section]
\newtheorem{defn}[thm]{Definition}
\newtheorem{corollary}[thm]{Corollary}

\newtheorem{lemma}[thm]{Lemma}

\newtheorem{remark}[thm]{Remark}
\newtheorem{example}[thm]{Example}
\newtheorem{assumption}[thm]{Assumption}

\newcommand{\pf}{\noindent{\bf Proof.} }
\def\qed{{\hfill $\Box$ \bigskip}}

\newcommand\cbrk{\text{$]$\kern-.15em$]$}}
\newcommand\opar{\text{\,\raise.2ex\hbox{${\scriptstyle
|}$}\kern-.34em$($}}
\newcommand\cpar{\text{$)$\kern-.34em\raise.2ex\hbox{${\scriptstyle |}$}}\,}

\def\<{\langle}
\def\>{\rangle}

\def\E{{\mathbb E}}

\newcommand\bL{\mathbb{L}}
\newcommand\bR{\mathbb{R}}
\newcommand\bH{\mathbb{H}}
\newcommand\bZ{\mathbb{Z}}

\newcommand\bD{\mathbb{D}}
\newcommand\bS{\mathbb{S}}

\newcommand\cF{\mathcal{F}}
\newcommand\cH{\mathcal{H}}

\newcommand\cP{\mathcal{P}}
\newcommand\cR{\mathcal{R}}

\newcommand\cO{\mathcal{O}}

\newcommand\frH{\mathfrak{H}}

\newcommand{\mysection}[1]{\section{#1}
\setcounter{equation}{0}}

\begin{document}

\title{\bf  A  weighted Sobolev space theory of  parabolic stochastic PDEs on non-smooth domains}

\author{Kyeong-Hun Kim}

\date{}

\author{Kyeong-Hun Kim\footnote{Department of
Mathematics, Korea University, 1 Anam-dong, Sungbuk-gu, Seoul, South
Korea 136-701, \,\, kyeonghun@korea.ac.kr. The research of this
author was supported by Basic Science Research Program through the
National Research Foundation of Korea(NRF) funded by the Ministry of
Education, Science and Technology (20110015961)}}

\maketitle

\begin{abstract}

In this paper we  study
 parabolic stochastic partial differential equations (see equation (\ref{5.15.0})) defined on arbitrary bounded domain $\cO \subset \bR^d$ allowing Hardy inequality:
  \begin{equation}
                          \label{eqn 5.2}
\int_{\cO}|\rho^{-1}g|^2\,dx\leq C\int_{\cO}|g_x|^2 dx, \quad
\forall g\in C^{\infty}_0(\cO),
\end{equation}
where $\rho(x)=\text{dist}(x,\partial \cO)$.
  Existence and uniqueness results are given in weighted
Sobolev spaces  $\frH^{\gamma}_{p,\theta}(\cO,T)$, where  $p\in [2,\infty)$, $\gamma\in \bR$ is the number of derivatives of solutions and $\theta$ controls the boundary behavior of solutions (see Definition \ref{def 1}). Furthermore  several
H\"older estimates of the solutions are also obtained. It is
allowed that the coefficients of the equations blow up near the
boundary.

\vspace*{.125in}

\noindent {\it Keywords: Hardy inequality, Stochastic  partial differential equation,
non-smooth domain, $L^p$-theory, weighted Sobolve space.}

\vspace*{.125in}

\noindent {\it AMS 2000 subject classifications:}  60H15, 35R60.

\end{abstract}

\maketitle

\mysection{Introduction}

  It is a classical result that
 Hardy inequality holds on Lipschitz domains (\cite{N}). There have been
 many other works concerning Hardy inequality. See e.g.  \cite{EH}, \cite{Wa} and
 references therein. We only mention that  inequality (\ref{eqn 5.2}) holds under much weaker condition than Lipschitz condition. For instance, it holds if
 $\cO$ has  plump complement, that is, there exist
 $b,\sigma\in (0,1]$ such that for any $s\in (0,\sigma]$ and $x\in \partial \cO$ there
 exists a point $y\in B_s(x) \cap \cO^c$ with
 $\text{dist}(y,\partial \cO)\geq bs$. For instance,
 $\cO_{\alpha}:=\{(x,y)\in \bR^2: x\in (-1,1), \,\,|x|^{\alpha}+|y|^{\alpha}<1\}$, where $\alpha\in (0,1)$, is a non-Lipschitz domain but   satisfies the plump complement condition.

Let $(\Omega,\cF,P)$ be a complete probability space,
$\{\cF_{t},t\geq0\}$ be an increasing filtration of $\sigma$-fields
$\cF_{t}\subset\cF$, each of which contains all $(\cF,P)$-null sets.
We assume that on $\Omega$ we are given independent  one-dimensional Wiener processes
 $w^{1}_{t},w^{2}_{t},...$
relative to $\{\cF_{t},t\geq0\}$.
The main goal of this article is to present an $L_p$-theory of  stochastic
 partial differential equation
\begin{equation}
                                          \label{5.15.0}
du=(a^{ij}u_{x^ix^j}+b^iu_{x^i}+
 cu+ f)\,dt
+(\sigma^{ik} u_{x^i}+\mu^{k}u+g^{k})\,dw^{k}_{t}
\end{equation}
given for $t>0$ and $x\in \cO$. Here
 $i$ and $j$ go from $1$ to $d$, and $k$ runs through
$\{1,2,...\}$ with the summation convention on $i,j,k$ being
enforced.
 The coefficients $a^{ij},b^i,c,\sigma^{ik},\mu^k$ and
the free terms $f,g^k$  are random functions depending on $t$ and
$x$. As mentioned in \cite{KL2}, such equations
 with a finite number of the processes $w^k_t$ appear, for instance, in
nonlinear filtering problems (estimations of the signal by observing
it when it is mixed with noises), and considering infinitely many
$w^k_t$ is instrumental in treating equations for measure-valued
processes, for instance, driven by space-time white noise (cf. \cite{Kr99}).

Equation (\ref{5.15.0})  has been extensively
studied by so many authors (see e.g. \cite{F, Kim04, Kim03, KK1, Kr00, Kr99, Kr96, KL2, Lo, MR, MP, PA, Yoo99} and references therein). We  give a very
brief review  only on the $L_p$-theory of the equation. The
$L_p$-theory ($p\geq 2$) of equation (\ref{5.15.0}) defined in
$\bR^d$ was introduced  by Krylov (\cite{Kr99}, \cite{Kr96}), and
later Krylov and Lototsky (\cite{KL2},\cite{KL1}) developed a weighted
$L_p$-theory of the equation defined on a
half space. It turned out that for SPDEs defined on domains the H\"older space approach does not allow
one to obtain results of reasonable generality, and the Sobolev
spaces without weights are trivially inappropriate.  Recently, these weighted $L_p$-theory  on half space were extended to    equations on smooth domains (e.g. \cite{Kim04, Kim03, KK1, KK2, Lo2}) and on (non-smooth) Lipschitz domain (\cite{Kim08}).

On non-smooth domains the spatial derivatives of the solution usually have additional singularities at the boundary which are due to the shape of the domain, see e.g.\ \cite{Gri85, Gri92} for the case of deterministic equations on polygonal domains and \cite{Lin} for a generalization to the stochastic setting. In the context of numerical approximation this suggests the use of non-uniform schemes. In \cite{Cio}   results of \cite{Kim08}  are used to prove that  the convergence rates of adequate non-uniform discretization schemes   are closely connected to the regularity of the solution measured in weighted Sobolev spaces.

 However, we acknowledge  that there is a   gap
  in the proof of Lemma 3.1 of \cite{Kim08}, and the main results of \cite{Kim08} are false unless  stronger assumption on the range of weights is assumed.
  We show this with a counterexample.
   In this article we reconstruct  the results in \cite{Kim08} under much weaker assumption on $\partial \cO$, but with smaller range of weights.
   The arguments used in this article are slightly different from those in \cite{Kim08}. For instance,  we do not use
any argument of flattening the boundary, which is a key tool  in  \cite{Kim08}. Most of our important steps are based just on the Hardy inequality and I\^o's formula.

As in \cite{Kim08, Kim04, Kim03, KK1, KK2, KL1, KL2,Lo2}  we prove the existence and
uniqueness results in weighted Sobolev classes $\frH^{\gamma}_{p,\theta}(\cO,T)$,  where $\gamma\in \bR$ is the number of derivatives of solutions and $\theta$ controls the boundary behavior of solutions (see Definition \ref{def 1}).  Also several
(interior) H\"older estimates of the solutions are also obtained (see Corollary \ref{corollary 1}).

As usual $\bR^{d}$ stands for the Euclidean space of points
$x=(x^{1},...,x^{d})$, $\bR^d_+=\{x\in \bR^d:x^1>0\}$ and
$B_r(x):=\{y\in \bR^d:|x-y|<r\}$.
 For $i=1,...,d$, multi-indices $\beta=(\beta_{1},...,\beta_{d})$,
$\beta_{i}\in\{0,1,2,...\}$, and functions $u(x)$ we set
$$
u_{x^{i}}=\partial u/\partial x^{i}=D_{i}u,\quad
D^{\beta}u=D_{1}^{\beta_{1}}\cdot...\cdot D^{\beta_{d}}_{d}u,
\quad|\beta|=\beta_{1}+...+\beta_{d}.
$$
We also use the notation $D^m$ for a partial derivative of order $m$
with respect to $x$. If we write $N=N(...)$, this means that the
constant $N$ depends only on what are in parenthesis. Throughout the
article, for functions depending on $\omega,t$ and $x$, the argument
$\omega \in \Omega$ will be omitted.

The author is grateful   to Ildoo Kim for  carefully reading  the earlier version of the  article and finding several typos and to N.V. Kryolv for providing the author an example.  The author is also thankful to
P.A Cioica and F. Lindner for useful discussions regarding the numerical approximations of SPDEs on non-smooth domains.

\mysection{Main results}
                                                      \label{section 03.30.03}

First we introduce some
  Sobolev spaces (see e.g
 \cite{Kr99}, \cite{kr99} and \cite{Lo2} for more details).
 Let $p\in (1,\infty)$,
$\gamma \in \bR$ and
$H^{\gamma}_p=H^{\gamma}_p(\bR^d)=(1-\Delta)^{-\gamma/2}L_p$ be the
set of all distributions $u$ such that $(1-\Delta)^{\gamma/2}u\in
L_p$. Define
$$
\|u\|_{H^{\gamma}_p}=\|(1-\Delta)^{\gamma/2}u\|_{L_p}:=\|\cF^{-1}[(1+|\xi|^2)^{\gamma/2}\cF(u)(\xi)]\|_p,
$$
where $\cF$ is the Fourier transform.
It is well known that if $\gamma$ is a nonnegative integer then
$$
H^{\gamma}_p=H^{\gamma}_{p}(\bR^d)=\{u:u,Du,...,D^{\gamma}u\in
L_p\}.
$$
Denote $\rho(x)=\text{dist}(x,\partial \cO)$ and fix    a bounded infinitely differentiable function $\psi$ defined
in $\cO$ such that (see e.g.  Lemma 4.13 in \cite{Ku} or  formula (2.6) in \cite{Lo3})
\begin{equation}
                                 \label{eqn 2.26.1}
N^{-1}\rho(x)\leq \psi(x)\leq N \rho(x),\quad \rho^{m}|D^{m}\psi_x|\leq
N(m)<\infty.
\end{equation}
Let $\zeta\in C^{\infty}_{0}(\bR_{+})$ be a   nonnegative function
satisfying
\begin{equation}
                                                       \label{11.4.1}
\sum_{n=-\infty}^{\infty}\zeta(e^{n+t})>c>0,\quad\forall t\in\bR.
\end{equation}
Note that any non-negative smooth function $\zeta\in C^{\infty}_0(\bR_+)$ so that  $\zeta>0$ on $[e^{-1},e]$ satisfies (\ref{11.4.1}).
For $x\in \cO$ and $n\in\bZ:=\{0,\pm1,...\}$ define
$$
\zeta_{n}(x)=\zeta(e^{n}\psi(x)).
$$
Then    $\text{supp}\,\zeta_n \subset \{x\in \cO:
e^{-n-k_0}<\rho(x)<e^{-n+k_0}\}=:G_n$ for some integer $k_0>0$,
\begin{equation}
                                   \label{eqn 07.29.2}
\sum_{n=-\infty}^{\infty}\zeta_{n}(x)\geq \delta >0,
\end{equation}
\begin{equation}
                                                         \label{10.10.06}
\zeta_{n} \in C^{\infty}_0(G_n), \quad |D^m \zeta_n(x)|\leq
N(\zeta,m) e^{mn}.
\end{equation}
For $p \geq 1$ and $\gamma \in \bR$, by $H^{\gamma}_{p,\theta}(\cO)$ we denote the set of all
distributions $u$ on $\cO$ such that
\begin{equation}
                                                 \label{10.10.03}
\|u\|_{H^{\gamma}_{p,\theta}(\cO)}^{p}:= \sum_{n\in\bZ}
e^{n\theta} \|\zeta_{-n}(e^{n} \cdot)u(e^{n}
\cdot)\|^p_{H^{\gamma}_p} < \infty.
\end{equation}
We also use the above notation for    $\ell_2$-valued functions
$g=(g_1,g_2,...)$, that is,
$$
\|g\|_{H^{\gamma}_p}=\|g\|_{H^{\gamma}_p(\ell_2)}=\||(1-\Delta)^{\gamma/2}g|_{\ell_2}\|_{L_p},
$$
$$
\|g\|^p_{H^{\gamma}_p(\cO,\ell_2)}=\sum_{n\in\bZ}
e^{n\theta} \|\zeta_{-n}(e^{n} \cdot)g(e^{n}
\cdot)\|^p_{H^{\gamma}_p(\ell_2)}.
$$
It is known (see Lemma \ref{lemma 4.16}) that if $\{\bar{\zeta}_n, n\in \bZ\}$ is another set of functions satisfying (\ref{eqn 07.29.2})
and (\ref{10.10.06}) (such functions can be easily constructed by mollifying the indicator functions $I_{G_n}$), then it yields the same space $H^{\gamma}_{p,\theta}(\cO)$.
Also  if $\gamma=n$ is a
nonnegative integer then
$$
L_{p,\theta}(\cO):=H^0_{p,\theta}(\cO)=L_p(\cO,\rho^{\theta-d}
dx),
$$
$$
H^{n}_{p,\theta}(\cO):=\{u:u,\rho Du,...,\rho^{n}D^{n}u\in
L_{p,\theta}(\cO)\},
$$
\begin{equation}
                        \label{eqn 2.14.1}
\|u\|^p_{H^{n}_{p,\theta}(\cO)}\sim \sum_{|\alpha|\leq
n}\int_{\cO}|\rho^{|\alpha|}D^{\alpha}u|^p \rho^{\theta-d}\,dx.
\end{equation}
We remark that  the space $H^{n}_{p,\theta}(\cO)$ is different
from $W^{n,p}(\cO, \rho,\varepsilon)$ introduced in \cite{Ku},
where
$$
W^{n,p}(\cO, \rho, \varepsilon)=\{u: u, Du, ...,D^nu \in
L_p(\cO, \rho^{\varepsilon}dx)\}.
$$

\noindent
Denote $\rho(x,y)=\rho(x)\wedge \rho(y)$. For
  $\nu \in(0,1]$
 and $k=0,1,2,...$, as in  \cite{GT}, define
$$
[f]^{(0)}_{k}=[f]^{(0)}_{k,\cO} =\sup_{\substack{x\in \cO\\
|\beta|=k}}\rho^{k}(x)|D^{\beta}f(x)|, \quad \quad
[f]^{(0)}_{k+\nu} =\sup_{\substack{x,y\in \cO\\ |\beta|=k}}
\rho^{k+\nu}(x,y)\frac{|D^{\beta}f(x)-D^{\beta}f(y)|}
{|x-y|^{\nu}},
$$
$$
|f|^{(0)}_{k}=\sum_{j=0}^{k}[f]^{(0)}_{j,\cO}, \quad \quad
|f|^{(0)}_{k+\nu}=|f|^{(0)}_{k}+ [f]^{(0)}_{k+\nu}.
$$
The above notation is used also for $\ell_2$ valued functions $g=(g^1,g^2,\cdots)$. For instance,
$$
[g]^{(0)}_{k} =\sup_{\substack{x\in \cO\\
|\beta|=k}}\rho^{k}(x)|D^{\beta}g(x)|_{\ell_2}.
$$

\noindent
   Here are some other  properties of the space
$H^{\gamma}_{p,\theta}(\cO)$ taken from
 \cite{Lo2} (also see \cite{Kr99-1}, \cite{kr99}).

\begin{lemma}
                                               \label{collection}

$(i)$ The space $C^{\infty}_0(\cO)$ is dense in
$H^{\gamma}_{p,\theta}(\cO)$.

$(ii)$ Assume that $\gamma-d/p=m+\nu$ for some $m=0,1,...$ and
$\nu\in (0,1]$, and $i,j$ are multi-indices such that $|i|\leq m,
|j|=m$. Then for any $u\in H^{\gamma}_{p,\theta}(\cO)$, we have
$$
\psi^{|i|+\theta/p}D^iu \in C(\cO), \quad
\psi^{m+\nu+\theta/p}D^ju\in C^{\nu}(\cO),
$$
$$
|\psi^{|i|+\theta/p}D^iu|_{C(\cO)}+[\psi^{m+\nu+\theta/p}D^ju]_{C^{\nu}(\cO)}\leq
C \|u\|_{ H^{\gamma}_{p,\theta}(\cO)}.
$$

$(iii)$ $\psi D, D \psi: H^{\gamma}_{p,\theta}(\cO)\to
H^{\gamma-1}_{p,\theta}(\cO)$ are bounded linear operators, and
for any $u\in H^{\gamma}_{p,\theta}(\cO)$
$$
\|u\|_{H^{\gamma}_{p,\theta}(\cO)} \leq N\|\psi
u_x\|_{H^{\gamma-1}_{p,\theta}(\cO)}+N\|u\|_{H^{\gamma-1}_{p,\theta}(\cO)}
\leq N \|u\|_{H^{\gamma}_{p,\theta}(\cO)},
$$
$$
\|u\|_{H^{\gamma}_{p,\theta}(\cO)} \leq N\|(\psi
u)_x\|_{H^{\gamma-1}_{p,\theta}(\cO)}+N\|u\|_{H^{\gamma-1}_{p,\theta}(\cO)}
\leq N \|u\|_{H^{\gamma}_{p,\theta}(\cO)}.
$$

$(iv)$ For any $\nu, \gamma \in \bR$,
$\psi^{\nu}H^{\gamma}_{p,\theta}(\cO)=H^{\gamma}_{p,\theta-p\nu}(\cO)$
and
\begin{equation}
                            \label{open}
\|u\|_{H^{\gamma}_{p,\theta-p\nu}(\cO)} \leq N
\|\psi^{-\nu}u\|_{H^{\gamma}_{p,\theta}(\cO)}\leq
N\|u\|_{H^{\gamma}_{p,\theta-p\nu}(\cO)}.
\end{equation}

$(v)$ If $\gamma\in (\gamma_0,\gamma_1)$ and $\theta\in
(\theta_0,\theta_1)$, then
$$
\|u\|_{H^{\gamma}_{p,\theta}(\cO)}\leq \varepsilon
\|u\|_{H^{\gamma_1}_{p,\theta}(\cO)}+
N(\gamma,p,\varepsilon)\|u\|_{H^{\gamma_0}_{p,\theta}(\cO)},
$$
$$
\|u\|_{H^{\gamma}_{p,\theta}(\cO)}\leq \varepsilon
\|u\|_{H^{\gamma}_{p,\theta_0}(\cO)}+
N(\gamma,p,\varepsilon)\|u\|_{H^{\gamma}_{p,\theta_1}(\cO)}.
$$

\end{lemma}

\begin{lemma}
                           \label{lemma 2.14}
(i)  Let $s=|\gamma|$ if $\gamma$ is an integer, and $s>|\gamma|$
otherwise, then
$$
\|a u\|_{H^{\gamma}_{p,\theta}(\cO)}\leq
N(d,s,\gamma)|a|^{(0)}_{s}\|u\|_{H^{\gamma}_{p,\theta}(\cO)}.
$$

 (ii) If $\gamma=0,1,2,...$, then
$$
\|a u\|_{H^{\gamma}_{p,\theta}(\cO)}\leq
N\sup_{\cO}|a|\|u\|_{H^{\gamma}_{p,\theta}(\cO )}
+N_0|a|^{(0)}_{\gamma}\|u\|_{H^{\gamma-1}_{p,\theta}(\cO)}
$$
where   $N_0=0$ if $\gamma=0$.

(iii) If $0\leq r\leq s$, then
$$
|a|^{(0)}_{r}\leq
N(d,r,s)(\sup_{\cO}|a|)^{1-r/s}(|a|^{(0)}_s)^{r/s}.
$$

The assertions also holds for $\ell_2$-valued functions $a$.

\end{lemma}

\pf
 For (i), see Theorem 3.1
in \cite{Lo2}. (ii) is an easy consequence of (\ref{eqn 2.14.1}),
and (iii) is from Proposition 4.2 in \cite{Li}.
\qed

\begin{remark}
                \label{remark 4.25}
By Lemma \ref{lemma 2.14}, for any $\nu\geq 0$, $\psi^{\nu}$ is a
point-wise multiplier in $H^{\gamma}_{p,\theta}(\cO)$. Thus if
$\theta_1\leq \theta_2$ then
$$
\|u\|_{H^{\gamma}_{p,\theta_2}(\cO)}\leq N
\|\psi^{(\theta_2-\theta_1)/p}u\|_{H^{\gamma}_{p,\theta_1}(\cO)}\leq
N\|u\|_{H^{\gamma}_{p,\theta_1}(\cO)}.
$$
\end{remark}

\begin{lemma}
                                       \label{lemma 4.16}
Let $\{\xi_n\}$ be a sequence of $C^{\infty}_0(\cO)$ functions
such that
$$|D^m\xi_n|\leq C(m) e^{nm}, \quad
\text{supp} \,\xi_n \subset \{x\in \cO:
e^{-n-k_0}<\rho(x)<e^{-n+k_0}\}
$$
 for some
$k_0>0$. Then for any $u\in H^{\gamma}_{p,\theta}(\cO)$,
$$
\sum_n e^{n\theta}\|\xi_{-n}(e^nx)u(e^nx)\|^p_{H^{\gamma}_{p}}\leq N
\|u\|^p_{H^{\gamma}_{p,\theta}(\cO)}.
$$
If in addition
$$
\sum_n |\xi_n|^p>\delta>0,
$$
then the reverse inequality also holds.

\end{lemma}

\pf
See Theorem 2.2 in \cite{Lo2}.
\qed

 Let $\cP$ be the predictable
$\sigma$-field generated by $\{\cF_{t},t\geq0\}$. Define
 $$
\bH^{\gamma}_{p}(T)=L_p(\Omega\times [0,T], \cP,H^{\gamma}_p), \quad \bH^{\gamma}_{p}(T,\ell_2)=L_p(\Omega\times [0,T], \cP,H^{\gamma}_p(\ell_2))
\quad
$$
$$
\bH^{\gamma}_{p,\theta}(\cO,T)= L_p(\Omega\times [0,T], \cP,H^{\gamma}_{p,\theta}(\cO)), \quad \bH^{\gamma}_{p,\theta}(\cO,T,\ell_2)= L_p(\Omega\times [0,T], \cP,H^{\gamma}_{p,\theta}(\cO,\ell_2)),
$$
$$
\bL_{p,\theta}(\cO,T)=\bH^{0}_{p,\theta}(\cO,T), \quad U^{\gamma}_p=L_p(\Omega,\cF_0,H^{\gamma-2/p}_p), \quad
U^{\gamma}_{p,\theta}(\cO)
=\psi^{-\frac{2}{p}+1}L_p(\Omega,\cF_0,H^{\gamma-2/p}_{p,\theta}(\cO)).
$$
That is, for instance, we say $u\in \bH^{\gamma}_{p,\theta}(\cO,T)$ if $u$ has a $H^{\gamma}_{p,\theta}(\cO)$-valued predictable version $\bar{u}$ defined on $\Omega \times [0,T]$ so that
$$
\|u\|_{\bH^{\gamma}_{p,\theta}(\cO,T)}=\|\bar{u}\|_{\bH^{\gamma}_{p,\theta}(\cO,T)}:=\left(\E \int^T_0 \|u(s,\cdot)\|^p_{H^{\gamma}_{p,\theta}(\cO)}dt\right)^{1/p}<\infty.
$$
Also by $u\in \psi^{-\frac{2}{p}+1}L_p(\Omega,\cF_0,H^{\gamma-2/p}_{p,\theta}(\cO))$ we mean $\psi^{2/p-1}u\in L_p(\Omega,\cF_0,H^{\gamma-2/p}_{p,\theta}(\cO))$, and
$$
\|u\|^p_{U^{\gamma}_{p,\theta}(\cO)}
:=\E\|\psi^{2/p-1}u\|^p_{H^{\gamma-2/p}_{p,\theta}(\cO)}.
$$

Below by $(u,\phi)$ we denote the image of $\phi\in C^{\infty}_0(\cO)$ under a distribution $u$.

\begin{defn}
            \label{def 1}
We write $u \in \frH^{\gamma+2}_{p,\theta}(\cO,T)$  if
 $u\in \psi \bH^{\gamma+2}_{p,\theta}(\cO,T)$,
$u(0,\cdot)\in U^{\gamma+2}_{p,\theta}(\cO)$ and
 for some
$f\in \psi^{-1}\bH^{\gamma}_{p,\theta}(\cO,T)$ and
 $g\in \bH^{\gamma+1}_{p,\theta}(\cO,T,\ell_2)$,
\begin{equation}
du=f\,dt +g^k \,dw^k_t,
\end{equation}
in the sense of distributions.  In other words, for any $\phi \in
C^{\infty}_{0}(\cO)$, the equality
$$
(u(t,\cdot),\phi)= (u(0,\cdot),\phi) + \int^{t}_{0}
(f(s,\cdot),\phi) \, ds + \sum^{\infty}_{k=1} \int^{t}_{0}
(g^k(s,\cdot),\phi)\, dw^k_s
$$
holds for all $t \leq T$ with probability $1$. In this situation
we  write $f=\bD u$ and $g=\bS u$.
The norm in $\frH^{\gamma+2}_{p,\theta}(\cO,T)$ is defined  by
$$
\|u\|_{\frH^{\gamma+2}_{p,\theta}(\cO,T)}=
\|\psi^{-1}u\|_{\bH^{\gamma+2}_{p,\theta}(\cO,T)} + \|\psi \bD
u\|_{\bH^{\gamma}_{p,\theta}(\cO,T)}
 +
\|\bS u\|_{\bH^{\gamma+1}_{p,\theta}(\cO,T,\ell_2)} +
\|u(0,\cdot)\|_{U^{\gamma+2}_{p,\theta}(\cO)}.
$$

\end{defn}

\begin{remark}

(i) Remember that for any $\alpha, \gamma \in \bR$,  $\|\psi^{\alpha} u\|_{H^{\gamma}_{p,\theta}(\cO)}\sim \|u\|_{H^{\gamma}_{p,\theta+p\alpha}(\cO)}$.
Thus the space
$\frH^{\gamma+2}_{p,\theta}(\cO,T)$ is independent of the choice of
$\psi$.

(ii) It is easy to check (see Remark 3.2 of \cite{Kr99} for details) that  for any $\phi\in C^{\infty}_0(\cO)$ and $g\in \bH^{\gamma+1}_{p,\theta}(\cO,T,\ell_2)$ we have $\sum_{k=1}^{\infty}\int^T_0 (g^k,\phi)^2ds <\infty$, and therefore the series of stochastic integral $\sum_{k=1}^{\infty}\int^t_0(g^k,\phi)dw^k_t$ converges in probability uniformly on $[0,T]$.
\end{remark}

\begin{thm}
             \label{thm cauchy}
Let $u_n \in \frH^{\gamma+2}_{p,\theta}(\cO,T), n=1,2,\cdots$ and $\|u_n\|_{\frH^{\gamma+2}_{p,\theta}(\cO,T)}\leq K$, where $K$ is a finite constant.
Then there exists a subsequence $n_k$ and a function $u\in \frH^{\gamma+2}_{p,\theta}(\cO,T)$ so that

(i) $u_{n_k}, u_{n_k}(0,\cdot), \bD u_{n_k}, \bS u_{n_k}$ converges weakly to $u, u(0,\cdot), \bD u$ and $\bS u$ in $\bH^{\gamma+2}_{p,\theta}(T,\cO), U^{\gamma+2}_{p,\theta}(\cO)$,  $\bH^{\gamma}_{p,\theta}(\cO)$  and  $\bH^{\gamma+1}_{p,\theta}(\cO,\ell_2)$ respectively;

(ii) for any $\phi\in C^{\infty}_0(\cO)$ and $t\in [0,T]$, we have $(u_{n_k}(t,\cdot),\phi) \to (u(t,\cdot),\phi)$ weakly in $L_p(\Omega)$.

\end{thm}

\pf
The proof is identical to that of the proof of Theorem 3.11 in \cite{Kr99}, where the theorem is proved when $\cO=\bR^d$.

\qed

\begin{thm}
               \label{lemma 2.6}
For any nonnegative integer $n\geq \gamma+2$, the set
$$
\frH^n_{p,\theta}(\cO,T)\bigcap  \bigcup^{\infty}_{k=1}
L_p(\Omega,C([0,T],C^n_0(\cO_k))),
$$
where $\cO_k:=\{x\in \cO: \psi(x)>1/k\}$, is  dense in
$\frH^{\gamma+2}_{p,\theta}(\cO,T)$.
\end{thm}

\pf
It is enough to repeat the proof of Theorem 2.9 in \cite{KL2}, where
the lemma is proved when  $\cO=\bR^d_+$.
\qed

\begin{thm}
                    \label{lemma 2.20.1}

(i) Let $2/p<\alpha<\beta\leq 1$ and $u\in
\frH^{\gamma+2}_{p,\theta}(\cO,T)$, then
$$
\E[\psi^{\beta-1}u]^p_{C^{\alpha/2-1/p}([0,T],H^{\gamma+2-\beta}_{p,\theta}(\cO))}\leq
N T^{(\beta-\alpha)p/2}\|u\|_{\frH^{\gamma+2}_{p,\theta}(\cO,T)},
$$
where $N$ is independent of $T$ and $u$.

 (ii) Let $p\in [2,\infty)$ and $T<\infty$, then
\begin{equation}
                                \label{eqn 4.23}
\E \sup_{t\leq
T}\|u(t)\|^p_{H^{\gamma+1}_{p,\theta}(\cO)}\leq
N\|u\|^p_{\frH^{\gamma+2}_{p,\theta}(\cO,T)},
\end{equation}
where $N=N(d,p,\gamma,\theta,\cO,T)$ is non-decreasing function of $T$. In particular, for any $t\leq T$,
\begin{equation}
                     \label{gronwall}
\|u\|^p_{\bH^{\gamma+1}_{p,\theta}(\cO,t)}\leq
N\int^t_0\|u\|^p_{\frH^{\gamma+2}_{p,\theta}(\cO,s)}\,ds.
\end{equation}
\end{thm}

\pf
The theorem is  proved in \cite{Kim08}  on Lipschitz domains, and the proof works on any arbitrary domains.

 (i). Due to the definition of $\E[\psi^{\beta-1}u]^p_{C^{\alpha/2-1/p}([0,T],H^{\gamma+2-\beta}_{p,\theta}(\cO))}$ we may assume $u(0)=0$. Let $\bD u=f$ and $\bS u=g$. By (\ref{10.10.03}) and Lemma
\ref{collection}(iv),
$$
I:=\E[\psi^{\beta-1}u]^p_{C^{\alpha/2-1/p}([0,T],H^{\gamma+2-\beta}_{p,\theta}(\cO))}
$$
\begin{equation}
                                            \label{eqn 3.24}
\leq N\sum_n e^{n(\theta+p(\beta-1))}\E
[u(t,e^nx)\zeta_{-n}(e^nx)]^p_{C^{\alpha/2-1/p}([0,T],H^{\gamma+2-\beta}_{p})}.
\end{equation}
Denote $T_0:=T^{(\beta-\alpha)p/2}$. By Corollary 4.12 in \cite{KR01}, there exists a constant $N>0$,
independent of $T$ and $u$, so that for any $a>0$,
$$
\E[u(t,e^nx)\zeta_{-n}(e^nx)]^p_{C^{\alpha/2-1/p}([0,T],H^{\gamma+2-\beta}_{p})}
\leq
NT_0a^{\beta-1}(a\|u(t,e^nx)\zeta_{-n}(e^nx)\|^p_{\bH^{\gamma+2}_{p}(T)}
$$
$$
+a^{-1}\|f(t,e^n)\zeta_{-n}(e^nx)\|^p_{\bH^{\gamma}_{p}(T)}+\|g(t,e^n)\zeta_{-n}(e^nx)\|^p_{\bH^{\gamma+1}_{p}(T,\ell_2)}).
$$
Take $a=e^{-np}$, then
(\ref{eqn 3.24}) yields
\begin{eqnarray*}
I &\leq& NT_0 (\sum_n
e^{n(\theta-p)}\|u(t,e^nx)\zeta_{-n}(e^nx)\|^p_{\bH^{\gamma+2}_p(T)}
+ \sum_n e^{n(\theta+p)}
\|f(t,e^nx)\zeta_{-n}(e^nx)\|^p_{\bH^{\gamma}_{p}(T)}\\
&+&\sum_n e^{n\theta}\|g(t,e^nx)\zeta_{-n}(e^nx)\|^p_{\bH^{\gamma+1}_p(T,\ell_2)})\\
&=& N T_0 \left(\|u\|^p_{\bH^{\gamma+2}_{p,\theta-p}(\cO,T)}
+\|f\|^p_{\bH^{\gamma}_{p,\theta+p}(\cO,T)} +\|g\|^p_{\bH^{\gamma+1}_{p,\theta}(\cO,T,\ell_2)}\right)
\leq N T_0\|u\|^p_{\frH^{\gamma+2}_{p,\theta}(\cO,T)}.
\end{eqnarray*}
Thus (i) is proved.

(ii). If $p>2$, (ii)  follows from (i). But for
the case $p=2$, we prove this differently. Obviously
$$
\E\sup_{t\leq
T}\|u(t)\|^p_{H^{\gamma+1}_{p,\theta}(\cO)}\leq N
\sum_ne^{n\theta} \E \sup_{t\leq T}
\|u(t,e^nx)\zeta_{-n}(e^nx)\|^p_{H^{\gamma+1}_{p}}.
$$
Note that $u_0\in U^{\gamma+2}_{p,\theta}(\cO)\subset L_p(\Omega,H^{\gamma+1}_{p,\theta}(\cO))$ since $p\geq 2$. By Remark 4.14 in \cite{KR01} with $\beta=1$ there, for any $a>0$,
$$
\E \sup_{t\leq T}
\|u(t,e^nx)\zeta_{-n}(e^nx)\|^p_{H^{\gamma+1}_{p}}\leq
N (a\|u(t,e^nx)\zeta_{-n}(e^nx)\|^p_{\bH^{\gamma+2}_{p}(T)}
$$
$$
+a^{-1}\|f(t,e^n)\zeta_{-n}(e^nx)\|^p_{\bH^{\gamma}_{p}(T)}+\|g(t,e^nx)\zeta_{-n}(e^nx)\|^p_{\bH^{\gamma+1}_{p}(T,\ell_2)}+
+\E\|u_0(e^nx)\zeta_{-n}(e^nx)\|^p_{H^{\gamma+1}_p}).
$$
Take $a=e^{-np}$ to get
\begin{eqnarray*}
\E\sup_{t\leq
T}\|u(t)\|^p_{H^{\gamma+1}_{p,\theta}(\cO)} &\leq&
N (\sum_n
e^{n(\theta-p)}\|u(t,e^nx)\zeta_{-n}(e^nx)\|^p_{\bH^{\gamma+2}_p(T)}
+ \sum_n e^{n(\theta+p)}
\|f(t,e^nx)\zeta_{-n}(e^nx)\|^p_{\bH^{\gamma}_{p}(T)}\\
&&+ \sum_n e^{n\theta}
\|g(t,e^nx)\zeta_{-n}(e^nx)\|^p_{\bH^{\gamma+1}_{p}(T)} +
\E\sum_ne^{n\theta}\|u_0(e^nx)\zeta_{-n}(e^nx)\|^p_{H^{\gamma+1}_{p}})\\
&=& N (\|u\|^p_{\bH^{\gamma+2}_{p,\theta-p}(\cO,T)}
+\|f\|^p_{\bH^{\gamma}_{p,\theta+p}(\cO,T)}+\|g\|^p_{\bH^{\gamma+1}_{p,\theta}(\cO)}+
\E\|u_0\|^p_{H^{\gamma+1}_{p,\theta}(\cO)})\\
&\leq& N\|u\|^p_{\frH^{\gamma+2}_{p,\theta}(\cO,T)}.
\end{eqnarray*}
Finally,
$$
\|u\|^p_{\bH^{\gamma+1}_{p,\theta}(\cO,t)}=\E \int^t_0 \|u\|^p_{H^{\gamma+1}_{p,\theta}(\cO)}ds\leq \int^t_0  (\E \sup_{r\leq s} \|u(r)\|^p_{H^{\gamma+1}_{p,\theta}(\cO)})ds\leq N \int^t_0 \|u\|^p_{\frH^{\gamma+2}_{p,\theta}(\cO,s)}ds.
$$
The theorem is proved.
\qed

\noindent
Fix a nonnegative constant $\varepsilon_0=\varepsilon(\gamma) \geq 0$ so that $\varepsilon_0>0$ only if $\gamma$ is not integer, and define
 $\gamma_+=|\gamma|$ if $\gamma$ is an integer, and
 $\gamma_+=|\gamma|+\varepsilon_0$ otherwise. Now we state our assumptions on the coefficients.

\begin{assumption}
                                \label{assumption 1}
(i) For each $x$, the coefficients $a^{ij}(t,x)$, $b^i(t,x)$
 $c(t,x)$, $\sigma^{ik}(t,x)$ and
 $\mu^k(t,x)$ are predictable functions of $(\omega,t)$.

 (ii) The coefficients $a^{ij}, \sigma^i$ are uniformly continuous in
 $x$, that is,
  for any $\varepsilon >0$ there exists $\delta=\delta(\varepsilon)>0$ such that
 $$
 |a^{ij}(t,x)-a^{ij}(t,y)|+|\sigma^i(t,x)-\sigma^i(t,y)|_{\ell_2}\leq \varepsilon
 $$
 for each $\omega,t$, whenever $x,y \in \cO$ and $|x-y|\leq \delta$.

(iii) There exist constant $\delta_0,K>0$ such that  for any
$\omega,t,x$ and $\lambda \in \bR^d$,
\begin{equation}
                                                           \label{1.23.01}
\delta_0 |\lambda|^2 \leq \bar{a}^{ij}(t,x) \lambda^{i} \lambda^{j}
 \leq K |\lambda|^2,
\end{equation}
 where $\bar{a}^{ij}=a^{ij}-\frac{1}{2}(\sigma^i,\sigma^j)_{\ell_2}$.

(iv) For any $\omega,t$

\begin{equation}
                   \label{eqn 2.20.1}
|a^{ij}(t,\cdot)|^{(0)}_{\gamma_+}+|\psi b^{i}(t,\cdot)|^{(0)}_{\gamma_+}+|\psi^2 c(t,\cdot)|^{(0)}_{\gamma_+}+ |\sigma^{i}(t,\cdot)|^{(0)}_{(\gamma+1)_+} +
|\psi \mu(t,\cdot)|^{(0)}_{(\gamma+1)_+} \leq K,
\end{equation}
and if $\gamma=0$, then for some $\varepsilon>0$,
\begin{equation}
                 \label{see this}
|\sigma^{i}(t,\cdot)|^{(0)}_{1+\varepsilon}+ |\psi \mu|^{(0)}_{1+\varepsilon}\leq K, \quad \quad \forall \omega,t.
\end{equation}

(v) There is a control on the behavior of $b^i,c$ and $\mu^k$ near
$\partial \cO$, namely,
\begin{equation}
                  \label{eqn 2.20.2}
\lim_{\substack{\rho(x)\to0}} \sup_{\omega,t}
\left(\rho(x)|b^i(t,x)|+\rho^2(x)|c(t,x)| +\rho(x)|\mu(t,x)|_{\ell_2} \right) =0.
\end{equation}
\end{assumption}

\begin{remark}
                         \label{remark coefficients}
Conditions (\ref{eqn 2.20.1}) and (\ref{eqn 2.20.2}) allow the
coefficients $b^i,c$ and $\nu$ to be unbounded and to blow up near
the boundary. In particular, (\ref{eqn 2.20.2}) is satisfied if for
some $\varepsilon, N>0$,
$$
|b^i(t,x)|+|\mu(t,x)|_{\ell_2}\leq N\rho^{-1+\varepsilon}(x),\quad
|c(t,x)|\leq N\rho^{-2+\varepsilon}(x).
$$

\end{remark}

The proof of following theorem is given in section \ref{section theorem 1}.

\begin{thm}
                                                     \label{theorem 1}
Let $p\in [2,\infty), \gamma\in [0,\infty), T<\infty$ and Assumption
 \ref{assumption 1} be satisfied. Then
there exists $\beta_0=\beta_0(p,d,\cO)>0$ so that if
\begin{equation}
                           \label{theta}
\theta\in
(p-2+d-\beta_0, p-2+d+\beta_0)
\end{equation}
 then  for any $f\in
\psi^{-1}\bH^{\gamma}_{p,\theta}(\cO,T)$, $g\in \bH^{\gamma+1}_{p,\theta}(\cO,T)$ and $u_0\in
U^{\gamma+2}_{p,\theta}(\cO)$ equation (\ref{5.15.0}) with
initial data $u_0$
 admits a unique solution $u$   in the class $\frH^{\gamma+2}_{p,\theta}(\cO,T)$,
and  for this solution
\begin{equation}
                                                          \label{eqn main}
\|u\|_{\frH^{\gamma+2}_{p,\theta}(\cO,T)}\leq C
(\|\psi f\|_{\bH^{\gamma}_{p,\theta}(\cO,T)}+\|g\|_{\bH^{\gamma+1}_{p,\theta}(\cO,T,\ell_2)}
+\|u_0\|_{U^{\gamma+2}_{p,\theta}(\cO)}),
\end{equation}
where  $C=C(d,p,\gamma,\theta,\delta_0,K,T,\cO)$.
\end{thm}

\begin{remark}
Note that Theorem \ref{theorem 1} is proved only for $\gamma\geq 0$. However the theorem can be extended for any $\gamma \in \bR$ by using  results for $\gamma \geq 0$ and arguments used e.g. in the proof of Theorem 2.16 of  \cite{KK1} (cf. \cite{KK2,KL2}). One difference is that, in place of Theorem 2.8 of \cite{kr99}, one has to use the corresponding version on bounded domains (Theorem 5.1 of \cite{Lo2}).
\end{remark}



 Lemma \ref{collection}(ii)  and Theorem \ref{lemma 2.20.1} easily  yield
  the following result.

\begin{corollary}
                          \label{corollary 1}
Let $u\in \frH^{\gamma+2}_{p,\theta}(\cO,\tau)$ be the solution in
Theorem \ref{theorem 1} (or in Theorem \ref{theorem 22} below).

(i) If $\gamma+2-d/p=m+\nu$ for some $m=0,1,...$, $\nu\in (0,1]$,
and $i,j$ are multi-indices such that $|i|\leq m, |j|=m$, then for
each $\omega,t$
$$
\psi^{|i|-1+\theta/p}D^iu\in C(\cO), \quad
\psi^{m-1+\nu+\theta/p}D^ju\in C^{\nu}(\cO).
$$
In particular,
$$
|\psi^{|i|}D^iu(x)|\leq N \psi^{1-\theta/p}(x).
$$

(ii) Let
$$
2/p<\alpha<\beta\leq 1, \quad \gamma+2-\beta-d/p=k+\varepsilon
$$
where $k\in \{0,1,2,...\}$ and $\varepsilon\in (0,1]$. Denote
$\delta=\beta-1+\theta/p$. Then for any multi-indices $i$ and $j$
such that $|i|\leq k$ and $|j|=k$, we have
$$
\E\sup_{t,s\leq
\tau}|t-s|^{-(p\alpha/2-1)}(|\psi^{\delta+|i|}D^i(u(t)-u(s))|^p_{C(\cO)}
$$
$$
+[\psi^{\delta+|j|+\varepsilon}D^j(u(t)-u(s))]^p_{C^{\varepsilon}(\cO)})<\infty.
$$

\end{corollary}

   Note that if $p=2$ then (\ref{theta}) is $\theta\in (d-\beta_0, d+\beta_0)$, but it is not clear whether  $d$ is included
    in the interval  in (\ref{theta}) if $p>2$  because  $\beta_0$ depends also on $p$.  Below we give  positive answer  if $p$ is close to $2$ and negative one for large $p$. For instance, {\bf{if
    $\theta=d=2$ then  in general Theorem \ref{theorem 1} do not hold for all $p>4$}}.

\vspace{3mm}

The proof of following theorem is given in section \ref{section theorem 22}.
\begin{thm}
                      \label{theorem 22}
 There exists $p_0>2$ so
that if $p\in [2,p_0)$ then there exists $\beta_1>0$ so that the assertion of Theorem \ref{theorem 1} holds for any  $\theta\in(d-\beta_1, d+\beta_1)$.
\end{thm}

\begin{remark}
Since $H^1_{p,d-p}(\cO)=\stackrel{\circ}{W^1_p}(\cO):=\{u: u,u_x\in L_p(\cO) \,\, \text{and}\,\, u|_{\partial\cO}=0\}$, Theorem \ref{theorem 22} with $\gamma \geq -1$ and $\theta \in (d-\beta_1,d]$ implies that there exists a
unique solution $u\in L_p(\Omega\times [0,T], \stackrel{\circ}{W^1_p}(\cO))$ for any $p\in [2,p_0)$.

\end{remark}

The following example  is due to N.V. Krylov and shows that  Theorem \ref{theorem 1} can not hold unless $\theta$ is sufficiently large and that in general Theorem \ref{theorem 22} is false for all large $p$.

\begin{example}
Let $\alpha\in (1/2, 1)$ and denote
$$
G_{\alpha}=\{z=x+iy: |\text{arg} \,z|< \frac{\pi}{2\alpha} \}, \quad \cO_{\alpha}=G_{\alpha}\cap \{z: |z|<10\},
$$
where $\text{arg} \,z$ is defined as a function taking values in so that $[\pi, -\pi)$. Define $v(z)=v(x,y)=\text{Re}\, z^{\alpha}=|z|^{\alpha}\cos \alpha \theta$, where $\tan \theta=y/x$. Then $\Delta v=0$ in $G_{\alpha}$ and $v=0$ on $\partial G_{\alpha}$.  We claim that   for some $N=N(\alpha)>1$,
$$
N^{-1}|z|^{\alpha-1}\leq |\rho^{-1}v|\leq N |z|^{\alpha-1}, \quad |Dv|+|\rho D^2v|\leq N|z|^{\alpha-1}.
$$
Since the second assertion is easy to check we prove the first one. If $|\text{arg}\,z|< \frac{\pi}{2\alpha}-\frac{\pi}{2}$ then $\rho(z)=|z|$, $|z|^{\alpha}\cos(\frac{\pi}{2}-\alpha \frac{\pi}{2})\leq |v|\leq |z|^{\alpha}$ and the claim is obvious.
Also if $\frac{\pi}{2\alpha}-\frac{\pi}{2}\leq |\text{arg}\,z|< \frac{\pi}{2\alpha}$, then $\rho(z)=|z||\sin (\frac{\pi}{2\alpha}-\theta)|$ and   $\cos \alpha \theta/ |\sin (\frac{\pi}{2\alpha}-\theta|$
is comparable to $1$ in $\{z: \frac{\pi}{2\alpha}-\frac{\pi}{2}\leq |\text{arg}\,z|< \frac{\pi}{2\alpha}\}$.

It follows that
$$
\int_{\cO_{\alpha}}\left(|\rho^{-1}v|^p+|Dv|+|\rho D^2v|\right)  \rho^{\theta-2}dx <\infty  \quad \Leftrightarrow \quad \theta>p(1-\alpha),
$$
and
$$
\int_{\cO_{\alpha}} (|\rho v_x|^p+|\rho v|^p) \rho^{\theta-2}dx<\infty, \quad \forall \,\, \theta>0.
$$
Now choose a smooth function $\xi \in C^{\infty}_0(B_2(0))$ so that $\xi=1$ on $B_1(0)$, and define $u(t,x,y):=t\xi(x,y)v(x,y)$. Then
\begin{equation}
          \label{ex}
du=(\Delta u +f)dt,
\end{equation}
 where $f:= t(-2\xi_{x^i}v_{x^i}-v\Delta \xi)+\xi v$. Above calculations show that $\rho f\in \bL_{p,\theta}(\cO_{\alpha},T)$ for any $\theta>0$ and that $u\in \frH^{2}_{p,p}(\cO_{\alpha},T)$.
 By Theorem \ref{theorem 1} we conclude that $u$ is  the unique solution of the above equation in $\frH^{2}_{p,p}(\cO_{\alpha},T)$.
It also follows that {\bf{the existence result of Theorem \ref{theorem 1} in  $\frH^2_{p,\theta}(\cO_{\alpha},T)$ fails}} whenever
$$
\theta\leq p(1-\alpha),
$$ because if there is any solution $w\in \frH^2_{p,\theta}(\cO_{\alpha},T)$ then $w\in \frH^{2}_{p,p}(\cO_{\alpha},T)$ and therefore due to the uniqueness result in $\frH^{2}_{p,p}(\cO_{\alpha},T)$, we get $u=w$. But this is not possible since $\|\rho^{-1}u\|_{\bL_{p,\theta}(\cO_{\alpha},T)}=\infty$.
In particular, if $\theta=d=2$ and $p>4$ we can choose $\alpha$ close to $1/2$ so that $2\leq p(1-\alpha)$, and consequently this leads to the fact that in general Theorem \ref{theorem 1} does not holds if $p>4$.
\end{example}

\mysection{A priori estimate}

In this section we develop  some estimations of solutions of equation (\ref{5.15.0}).
First, we  introduce a result on SPDEs defined on entire space $\bR^d$.

\begin{lemma}
                                                    \label{lemma 3}
Let $a^{ij}$ and $\sigma^{ij}$ be independent of $x$.
 Also suppose that  $f\in \bH^{\gamma}_p(T)$, $g\in \bH^{\gamma+1}_{p}(T,\ell_2)$,
  $u_0\in U^{\gamma+2}_p$
 and
 $u\in \bH^{\gamma+1}_p(T)$ is a solution of
 \begin{equation}
                            \label{eqn 08.8}
 du=(a^{ij}u_{x^ix^j}+f)+ (\sigma^{ik}u_{x^i}+g^k)dw^k_t \quad
 u(0,\cdot)=u_0.
 \end{equation}
 Then  $u\in \bH^{\gamma+2}_p(T)$,  and
\begin{equation}
                                                   \label{eqn 3}
\|u\|^p_{\bH^{\gamma+2}_p(T)} \leq N (\|u\|^p_{\bH^{\gamma+1}_p(T)}
 +\|f\|^p_{\bH^{\gamma}_p(T)}+\|g\|^p_{\bH^{\gamma+1}_p(T,\ell_2)}
 +\|u_0\|^p_{U^{\gamma+2}_p}),
\end{equation}
where $N$ depends only on $d,p,\delta_0,K$ (not on $T$).

\end{lemma}
\pf
This is a well known result. By Theorem 4.10 in \cite{Kr99},
$$
\|u_{xx}\|^p_{\bH^{\gamma}_p(T)} \leq C(d,p) (\|f\|^p_{\bH^{\gamma}_p(T)}+\|g\|^p_{\bH^{\gamma+1}_p(T,\ell_2)}
 +\|u_0\|^p_{U^{\gamma+2}_p}).
 $$
 This and the relation $\|u\|_{H^{\gamma+2}_p}=\|(1-\Delta)u\|_{H^{\gamma}_p}\leq (\|u\|_{H^{\gamma}_p}+\|u_{xx}\|_{H^{\gamma}_p})$
certainly prove (\ref{eqn 3}).
\qed

In the following lemma there is no restriction on $\theta,\gamma$ and $\partial \cO$, that is $\theta, \gamma \in \bR$ and $\cO$ is any arbitrary domain.

\begin{lemma}
                            \label{lemma 4.2}
Let $a^{ij}$ and $\sigma^{ik}$ be independent of $x$. Suppose  $f\in
\psi^{-1}\bH^{\gamma}_{p,\theta}(\cO,T)$,  $g\in \bH^{\gamma+1}_{p,\theta}(\cO,T,\ell_2)$, $u_0\in
U^{\gamma+2}_{p,\theta}(\cO)$ and $u\in
\frH^{\gamma+1}_{p,\theta}(\cO,T)$ is a solution of the equation
$$
du=(a^{ij}u_{x^ix^j}+f)dt+(\sigma^{ik}u_{x^i}+g^k)dw^k_t, \quad \quad u(0,\cdot)=u_0.
$$
Then $u\in
\frH^{\gamma+2}_{p,\theta}(\Omega,T)$, and
\begin{equation}
                         \label{eqn 2.16.5}
\|\psi^{-1}u\|^p_{\bH^{\gamma+2}_{p,\theta}(\cO,T)}\leq
N(\|\psi^{-1}u\|^p_{\bH^{\gamma+1}_{p,\theta}(\cO,T)}+\|\psi
f\|^p_{\bH^{\gamma}_{p,\theta}(\cO,T)}+\|g\|^p_{\bH^{\gamma+1}_{p,\theta}(\cO,T,\ell_2)}
+\|u_0\|^p_{U^{\gamma+2}_{p,\theta}(\cO)}).
\end{equation}
\end{lemma}

\pf We just repeat the arguments used in \cite{Kim08} on Lipschitz domains.
Remember that  by Lemma \ref{collection} we have $\|\psi^{-1}u\|_{H^{\gamma+2}_{p,\theta}(\cO)}\sim \|u\|_{H^{\gamma+2}_{p,\theta-p}(\cO)}$. Thus,
$$
\|\psi^{-1}u\|^p_{\bH^{\gamma+2}_{p,\theta}(\cO,T)}\leq N
\sum_{n=-\infty}^{\infty}e^{n(\theta-p)}\|u(t,e^nx)\zeta_{-n}(e^nx)\|^p_{\bH^{\gamma+2}_p(T)}
$$
\begin{equation}
                    \label{eqn 2.16.3}
=N\sum_{n=-\infty}^{\infty}e^{n(\theta-p+2)}\|v_n\|^p_{\bH^{\gamma+2}_p(e^{-2n}T)},
\end{equation}
where $v_n(t,x):=u(e^{2n}t,e^nx)\zeta_{-n}(e^nx)$.
 Note that since $v_n$ has compact support in $\bR^d$ and can be regarded as distribution defined on $\bR^d$. Thus we conclude $v_n\in
\bH^{\gamma+1}_p(e^{-2n}T)$. Also note that it  satisfies
$$
dv_n=(a^{ij}(e^{2n}t)v_{nx^ix^j}+f_n)dt+(\sigma^{ik}(e^{2n}t)v_{nx^i}+g^k_n)dw^k(n)_t, \quad v_n(0)=u_0(e^nx)\zeta_{-n}(e^nx),
$$
where $w^k(n)_t:=e^{-n}w^k_{e^nt}$ are independent Wiener processes,
\begin{eqnarray*}
f_n(t,x)&=&-2e^{n}a^{ij}(e^{2n}t,x)u_{x^i}(e^{2n}t,e^nx)e^n\zeta_{-nx^j}(e^nx)
 -a^{ij}u(e^{2n}t,e^nx)e^{2n}\zeta_{-nx^ix^j}(e^nx)\\
&+&e^{2n}f(e^{2n}t,e^nx)\zeta_{-n}(e^nx),
\end{eqnarray*}
and
$$
g^k_n=-\sigma^{ik}(e^{2n}t)u(e^{2n}t,e^nx)e^n\zeta_{-nx^i}(e^nx)+e^ng^k(e^{2n}t,e^nx)\zeta_{-n}(e^nx).
$$
Since $\zeta_{-n}$ has compact support in $\cO$ and $u\in \bH^{\gamma+1}_{p,\theta}(\cO,T)$, we easily check that
 $$
 f_n\in \bH^{\gamma}_p(e^{-2n}T), \quad g_n\in \bH^{\gamma+1}_p(e^{-2n}T,\ell_2).
 $$
 Thus by
Lemma \ref{lemma 3}, we have $v_n \in \bH^{\gamma+2}_p(e^{2n}T)$
and
$$
\|v_{n}\|^p_{\bH^{\gamma+2}_p(e^{2n}T)}\leq
N(\|v_{n}\|^p_{\bH^{\gamma+1}_p(e^{2n}T)}+\|f_n\|^p_{\bH^{\gamma}_p(e^{2n}T)}+\|g_n\|^p_{\bH^{\gamma+1}_{p}(e^{2n}T,\ell_2)}
+\|v_n(0)\|^p_{U^{\gamma+2}_p}),
$$
where $N=N(d,p,\gamma,\delta_0,K)$ is independent of $n$ and $T$. Next we apply Lemma
\ref{lemma 4.16} with $\xi_n=e^{-n}\zeta_{nx^i}$ or
$\xi_n=e^{-2n}\zeta_{nx^ix^j}$ and get
\begin{eqnarray*}
\sum_{n=-\infty}^{\infty}e^{n(\theta-p+2)}\|f_n\|^p_{\bH^{\gamma}_p(e^{-2n}T)}
&\leq& N \sum_n
e^{n\theta}\|u_{x^i}(t,e^nx)e^n\zeta_{-nx^j}(e^nx)\|^p_{\bH^{\gamma}_p(T)}\\
&+&N\sum_n
e^{n(\theta-p)}\|u(t,e^nx)e^{2n}\zeta_{-nx^ix^j}(e^nx)\|^p_{\bH^{\gamma}_p(T)}\\
 &+&N \sum_n
e^{n(\theta+p)}\|f(t,e^nx)\zeta_{-n}(e^nx)\|^p_{\bH^{\gamma}_p(T)}\\
&\leq& N\|u_x\|^p_{\bH^{\gamma}_{p,\theta}(\cO,T)}+
N\|u\|^p_{\bH^{\gamma}_{p,\theta-p}(\cO,T)} +N\|
f\|^p_{\bH^{\gamma}_{p,\theta+p}(\cO,T)}\\
&\leq& N\|\psi^{-1}u\|^p_{\bH^{\gamma+1}_{p,\theta}(\cO,T)}+N\|\psi f\|^p_{\bH^{\gamma}_{p,\theta}(\cO,T)}.
\end{eqnarray*}
Similarly,
\begin{eqnarray*}
&&\sum_{n=-\infty}^{\infty}e^{n(\theta-p+2)}\|g_n\|^p_{\bH^{\gamma+1}_p(e^{-2n}T)} \\
&\leq& N \sum_n
e^{n(\theta-p)}\|u(t,e^nx)e^{n}\zeta_{-nx}(e^nx)\|^p_{\bH^{\gamma+1}_p(T)}+N\sum_n
e^{n\theta}\|g(t,e^nx)\zeta_{-n}(e^nx)\|^p_{\bH^{\gamma+1}_p(T,\ell_2)}\\
&\leq& N\|\psi^{-1}u\|^p_{\bH^{\gamma+1}_{p,\theta}(\cO,T)}+N\|g\|^p_{\bH^{\gamma+1}_{p,\theta}(\cO,T,\ell_2)}.
\end{eqnarray*}
Also,
$$
\sum_{n=-\infty}^{\infty}e^{n(\theta-p+2)}\|v_n(0)\|^p_{U^{\gamma+2}_p}
\leq N\|u_0\|^p_{U^{\gamma+2}_{p,\theta}(\cO)}.
$$
Thus the lemma is proved.

\qed

\begin{remark}
                  \label{remark 08.5}
Let $\gamma \geq 0$. By (\ref{eqn 2.16.5}) and the inequality (see Lemma \ref{collection}(v))
$$
\|\psi^{-1}u\|_{H^{\gamma+1}_{p,\theta}(\cO)}\leq \varepsilon
\|\psi^{-1}u\|_{H^{\gamma+2}_{p,\theta}(\cO)}+N(\varepsilon)
\|\psi^{-1}u\|_{L_{p,\theta}(\cO)},
$$
we easily get
\begin{equation}
                         \label{eqn 07.29.1}
\|\psi^{-1}u\|^p_{\bH^{\gamma+2}_{p,\theta}(\cO,T)}\leq
N(\|\psi^{-1}u\|^p_{\bL_{p,\theta}(\cO,T)}+\|\psi
f\|^p_{\bH^{\gamma}_{p,\theta}(\cO,T)}+\|g\|^p_{\bH^{\gamma+1}_{p,\theta}(\cO,T,\ell_2)}
+\|u_0\|^p_{U^{\gamma+2}_{p,\theta}(\cO)}).
\end{equation}
This shows that to estimate $\|u\|_{\frH^{\gamma+2}_{p,\theta}(\cO)}$ it is enough to estimate  $\|\psi^{-1}u\|^p_{\bL_{p,\theta}(\cO,T)}$.
\end{remark}

In the following lemma we estimate $\|\psi^{-1}u\|^p_{\bL_{p,\theta}(\cO,T)}$ when $\theta=d-2+p$ using the Hardy inequality.

\begin{lemma}
                 \label{lemma 4.5}
Let $a^{ij}$ and $\sigma^{ik}$ be independent of $x$.  Then for any $u\in \frH^2_{p,d-2+p}(\cO,T)$,
we have
\begin{eqnarray}
                              \|u\|_{\frH^2_{p,d-2+p}(\cO,T)} &\leq& N
 \|\psi (\bD u-a^{ij}u_{x^ix^j})\|_{\bL_{p,d-2+p}(\cO,T)} \nonumber\\
 &+& N\|\bS u-\sigma^{i\cdot}u_{x^i}\|_{\bH^{1}_{p,d-2+p}(\cO,T,\ell_2)}
 +N\|u(0)\|_{U^2_{p,d-2+p}(\cO)}, \label{eqn 4.14.5}
 \end{eqnarray}
 where $N=N(d, p,\cO)$.
 \end{lemma}

 \pf
{\bf{Step 1}}. First assume that $u\in L_p(\Omega,C([0,T],C^2_0(\cO_k)))$ for some $k$,
where $\cO_k:=\{x\in \cO: \psi(x)>1/k\}$, so that $u$ is sufficiently smooth in $x$ and vanishes near the boundary $\partial \cO$. Denote
 $$
 f=\bD u-a^{ij}u_{x^ix^j}, \quad g=\bS u-\sigma^{ik}u_{x^i}, \quad u_0=u(0).
 $$
 Then for each $x\in \cO$,
 $$
 u(t,x)=u_0(x)+\int^t_0 (a^{ij}u_{x^ix^j}+f)ds+\int^t_0 (\sigma^{ik}u_{x^i}+g^k)dw^k_t,
 $$
 for all $t\leq T$ (a.s.).  Applying It\^o's formula to $|u(t,x)|^p$,
 \begin{eqnarray*}
 |u(T)|^p&=&|u_0|^p+ p\int^T_0
 |u|^{p-2}u(a^{ij}u_{x^ix^j}+f)\,dt +\int^T_0p|u|^{p-2}u(\sigma^{ik}u_{x^i}+g^k)dw^k_t\\
 &+&\frac{1}{2}p(p-1)\int^T_0|u|^{p-2}\sum_{k=1}^{\infty}(\sigma^{ik}u_{x^i}+g^k)^2dt.
 \end{eqnarray*}
 Note that
 $$
 \frac{1}{2}p(p-1)|u|^{p-2} \sum_{k=1}^{\infty}(\sigma^{ik}u_{x^i}+g^k)^2=p(p-1)|u|^{p-2}\left(\alpha^{ij}u_{x^i}u_{x^j}+u_{x^i}(\sigma^{i},g)_{\ell_2}+\frac{1}{2}|g|^2_{\ell_2}\right),
 $$
 where $\alpha^{ij}=\frac{1}{2}(\sigma^i,\sigma^j)_{\ell_2}$.
 Taking expectation, integrating over $\cO$ and doing integration by parts (that is, $\int_{\cO}p|u|^{p-2}ua^{ij}u_{x^ix^j}dx=-p(p-1)\int_{\cO}a^{ij}|u|^{p-2}u_{x^i}u_{x^j}dx$), we get
 \begin{eqnarray*}
 p(p-1)\E \int^T_0 \int_{\cO}\bar{a}^{ij}|u|^{p-2}u_{x^i}u_{x^j}dxdt &\leq& \E\int_{\cO}|u_0|^pdx +p\E\int^T_0\int_{\cO} |u|^{p-1}|f|dxdt \\
 &+& p(p-1)\int^t_0\int_{\cO}|u|^{p-2}(u_{x^i}(\sigma^{i},g)_{\ell_2}+\frac{1}{2}|g|^2_{\ell_2}) dxdt.
 \end{eqnarray*}
 Note that for each $\omega,t$ we have $v:=|u|^{p/2}\in \{f: f,f_x\in L_2(\cO), f|_{\partial \cO}=0\}$, and
 $v_x=\frac{p}{2}|u|^{p/2-2}uu_x$. Thus by Hardy Inequality (see (\ref{eqn 5.2})),
 \begin{equation}
                                    \label{eqn 08.2.1}
  \int_{\cO}|\psi^{-1}u|^{p}\psi^{p-2}dx=\int_{\cO}|\psi^{-1}v|^2 dx\leq N\int_{\cO}|v_x|^2 dx \leq N\int_{\cO}|u|^{p-2}|u_x|^2 dx.
  \end{equation}
 Also note that
 \begin{eqnarray*}
 \int_{\cO}|u|^{p-1}|f|dx=\int_{\cO}|\psi^{-1}u|^{p-1}|\psi f| \psi^{p-2}dx &\leq& \varepsilon \int_{\cO}|\psi^{-1}u|^p \psi^{p-2}dx +N(\varepsilon)\int_{\cO}|\psi f|^p \psi ^{p-2}dx,
 \end{eqnarray*}
 \begin{eqnarray*}
 \int_{\cO}|u|^{p-2}u_{x^i}(\sigma^{i},g)_{\ell_2}dx &\leq & N|\sigma|_{\ell_2}\int_{\cO}|u|^{p-2}|u_x||g|_{\ell_2}dx\\
 &=&N|\sigma|_{\ell_2}\int_{\cO}|\psi^{-1}u|^{p-2}|u_x||g|_{\ell_2} \psi^{p-2}dx\\
 &\leq& \varepsilon \int_{\cO}|\psi^{-1}u|^p\psi^{p-2}dx+\varepsilon \int_{\cO}|u_x|^p\psi^{p-2}dx+N(\varepsilon)\int_{\cO}|g|^p_{\ell_2} \psi^{p-2}dx.
 \end{eqnarray*}
 Similarly,
 $$
 \int_{\cO}|u|^{p-2}|g|^2_{\ell_2} dx \leq \varepsilon \int_{\cO}|\psi^{-1}u|^p\psi^{p-2}dx+N(\varepsilon)\int_{\cO}|g|^p_{\ell_2} \psi^{p-2}dx.
 $$
 Since $(\bar{a}^{ij})\geq \delta_0 I$, we have $\delta |u|^{p-2}|Du|^2\leq \bar{a}^{ij}|u|^{p-2}u_{x^i}u_{x^j}$, and therefore from above calculations
\begin{eqnarray*}
(1-N_0 \varepsilon )\E \int^T_0 \int_{\cO}|\psi^{-1}u|^{p}\psi^{p-2}dxdt
&\leq& N \E \int_{\cO}|\psi^{\frac{2}{p}-1}u(0)|^p
\psi^{p-2}\,dx
+N\varepsilon \E\int^T_0\int_{\cO}|u_x|^{p}\psi^{p-2}\,dx\,dt\\
&+&N(\varepsilon)\E\int^T_0\int_{\cO}|\psi f|^p\psi^{p-2}dxdt+ N(\varepsilon)\E\int^T_0\int_{\cO}|g|^p_{\ell_2}\psi^{p-2}dxdt.
\end{eqnarray*}
Thus for any $\varepsilon>0$ so that $\varepsilon N_0 <1/2$, we have
\begin{eqnarray}
\|\psi^{-1}u\|_{\bL_{p,d-2+p}(\cO,T)}&\leq& N\|u_0\|_{U^1_{p,d-2+p}}+ N\varepsilon \|u_x\|_{\bL_{p,d-2+p}(\cO,T)} \nonumber\\
&+&N(\varepsilon)\|\psi f\|_{\bL_{p,d-2+p}(\cO,T)}+N(\varepsilon)\|g\|_{\bL_{p,d-2+p}(\cO,T,\ell_2)}.                       \label{eqn 8.4.4}
         \end{eqnarray}
   This and (\ref{eqn 07.29.1}) easily lead to     (\ref{eqn 4.14.5}).

 {\bf{Step 2}}. General case. We use Theorem \ref{lemma 2.6}. Take a sequence $u^n\in \frH^{2}_{p,d-2+p}(\cO,T)$ so that
 $u^n  \to u$ in $\frH^2_{p,d-2+p}(\cO,T)$ and each $u^n\in   L_p(\Omega,C([0,T],C^2_0(G_k)))$ for some $k=k(n)$.  By Step 1, we have (\ref{eqn 4.14.5}) with $u^n$
 in place of $u$. Now it is enough to let  $n \to \infty$.
\qed

The following lemma virtually says that if Theorem \ref{theorem 1} holds for some $\theta_0\in \bR$, then it also holds for all $\theta$ near $\theta_0$.

\begin{lemma}
                                       \label{lemma 4.25}
Suppose that there exists a $\theta_0\in \bR$ so that
for any $u\in \frH^2_{p,\theta_0}(\cO,T)$ we have
\begin{equation}
                            \label{eqn 8.5.1}
\|u\|_{\frH^2_{p,\theta_0}(\cO,T)}\leq N\left(\|\psi \bD u-\psi a^{ij}u_{x^ix^j}\|_{\bL_{p,\theta_0}(\cO,T)}+\|\bS u-\sigma^iu_{x^i}\|_{\bH^1_{p,\theta_0}(\cO,T,\ell_2)}+\|u(0)\|_{U^2_{p,\theta_0}(\cO)}\right).
\end{equation}
Then there exists $\varepsilon_0=\varepsilon_0(N,\theta_0,p)>0$ so that for any $\theta\in (\theta_0-\varepsilon_0, \theta_0+\varepsilon_0)$ and
$v\in \frH^{2}_{p,\theta}(\cO,T)$ it holds that
$$
\|v\|_{\frH^2_{p,\theta}(\cO,T)}\leq N\left(\|\psi \bD v-\psi a^{ij}v_{x^ix^j}\|_{\bL_{p,\theta}(\cO,T)}+\|\bS v-\sigma^iv_{x^i}\|_{\bH^1_{p,\theta}(\cO,T,\ell_2)}+\|v(0)\|_{U^2_{p,\theta}(\cO)}\right).
$$
\end{lemma}

\pf  Let $v\in \frH^{2}_{p,\theta}(\cO,T)$.
Denote $\nu=(\theta_0-\theta)/{p}$ and  $u:=\psi^{\nu}v$, then by (\ref{open}), $u \in \frH^2_{p,\theta_0}(\cO,T)$. Also it is easy to check that $\bD u=\psi^{\nu}\bD v$, $\bS u=\psi^{\nu}\bS v$ and
$$
\bD u-a^{ij}u_{x^ix^j}=\psi^{\nu}(\bD v-a^{ij}v_{x^ix^j})-2a^{ij}v_{x^i}(\psi^{\nu})_{x^j}-a^{ij}v(\psi^{\nu})_{x^ix^j},
$$
$$
\bS u-\sigma^{i}u_{x^i}=\psi^{\nu}(\bS v-\sigma^iv_{x^i})-\sigma^iv (\psi^{\nu})_{x^i}.
$$
Note, since $\psi_x$ and $\psi \psi_{xx}$ are bounded, if $\nu\leq 1$ then
\begin{equation}
                       \label{eqn little}
|(\psi^{\nu})_{x^j}|=\nu |\psi^{\nu-1}\psi_{x^i}|\leq N\nu \psi^{\nu-1},\quad \quad
|(\psi^{\nu})_{x^ix^j}|\leq N\nu \psi^{\nu-2}.
\end{equation}
By assumption (see (\ref{eqn 8.5.1})) and (\ref{eqn little}))
\begin{eqnarray*}
\|\psi^{\nu}v\|_{\frH^{2}_{p,\theta_0}(\cO,T)}&\leq&
 N \|\psi^{\nu} \psi(\bD v-a^{ij}v_{x^ix^j})\|_{\bL_{p,\theta_0}(\cO,T)}+N\|\psi^{\nu}(\bS v-\sigma^iv_{x^i})\|_{\bH^1_{p,\theta_0}(\cO,T,\ell_2)}\\
 &+&N\nu (\|\psi^{\nu}\psi^{-1}v\|_{\bH^1_{p,\theta_0}(\cO,T)}+\|\psi^{\nu}v_{x}\|_{\bL_{p,\theta_0}(\cO,T)})
 +N\|\psi^{\nu}v(0)\|_{U^2_{p,\theta_0}(\cO)}.
\end{eqnarray*}
This certainly implies (see (\ref{open}))
\begin{eqnarray*}
\|v\|_{\frH^{2}_{p,\theta}(\cO,T)}&\leq&
 N \| \psi(\bD v-a^{ij}v_{x^ix^j})\|_{\bL_{p,\theta}(\cO,T)}+N\|\bS v-\sigma^iv_{x^i}\|_{\bH^1_{p,\theta}(\cO,T,\ell_2)}\\
 &+&N_1\nu \left(\|\psi^{-1}v\|_{\bH^1_{p,\theta}(\cO,T)}+\|v_{x}\|_{\bL_{p,\theta}(\cO,T)}\right)+N \|v(0)\|_{U^2_{p,\theta}(\cO)}.
\end{eqnarray*}
It follows that  the claim of the lemma holds  for all sufficiently small $\nu$, that is for any $\theta$ so that $N_1|\theta_0-\theta|/p <1$.
The lemma is proved.
\qed

Remark \ref{remark 08.5}, Lemma \ref{lemma 4.5} and Lemma \ref{lemma 4.25}  obviously lead to the following result.

\begin{corollary}
                     \label{corollary good}
Suppose that  $\gamma\geq 0$ and the coefficients $a^{ij}, \sigma^{ik}$ are  independent of $x$.
Then there exists $\beta_0=\beta_0(d,p,\cO)>0$ so that  if $\theta\in (d-2+p-\beta_0, d-2+p+\beta_0)$, $f\in \bH^{\gamma}_{p,\theta}(\cO,T)$, $g\in \bH^{\gamma+1}_{p,\theta}(\cO,T,\ell_2)$, $u_0\in U^{\gamma+2}_{p,\theta}(\cO)$ and  $u\in \frH^{\gamma+2}_{p,\theta}(\cO,T)$ is a solution of (\ref{eqn 08.8}), then we have
\begin{equation}
                  \label{eqn 08.9.1}
\|u\|_{\frH^{\gamma+2}_{p,\theta}(\cO,T)}\leq N\left(\|\psi f\|_{\bH^{\gamma}_{p,\theta}(\cO,T)}+\|g\|_{\bH^{\gamma+1}_{p,\theta}(\cO,T,\ell_2)}+\|u_0\|_{U^{\gamma+2}_{p,\theta}(\cO)}\right),
\end{equation}
where $N=N(d,p,\theta,\delta_0,K,\cO,T)$.
\end{corollary}

Now we prove a priori estimate for solutions of the equation
\begin{equation}
                         \label{eqn general}
du=(a^{ij}u_{x^ix^j}+b^iu_{x^i}+
 cu+ f)\,dt
+(\sigma^{ik} u_{x^i}+\mu^{k}u+g^{k})\,dw^{k}_{t}, \quad u(0)=u_0.
\end{equation}

\begin{thm}
                  \label{a priori}
Suppose $\gamma \geq 0$, $\theta\in (d-2+p-\beta_0, d-2+p+\beta_0)$ and Assumption \ref{assumption 1} are satisfied. Also let $f\in \bH^{\gamma}_{p,\theta}(\cO,T)$, $g\in \bH^{\gamma+1}_{p,\theta}(\cO,T,\ell_2)$ and  $u_0\in U^{\gamma+2}_{p,\theta}(\cO)$. Then
estimate (\ref{eqn 08.9.1}) holds given that $u\in \frH^{\gamma+2}_{p,\theta}(\cO,T)$ is a solution of (\ref{eqn general}).
\end{thm}

\pf  {\bf{Step 1}}. Assume
$$
|a^{ij}(t,x)-a^{ij}(t,y)|+|\sigma^i(t,x)-\sigma^{i}(t,y)|_{\ell_2}+|\psi(x)b^i(t,x)|+|\psi^2(x)c(t,x)|+|\psi \mu|_{\ell_2}\leq \kappa, \quad
\forall \omega,t,x,y.
$$
We prove that there exists
$\kappa_0=\kappa_0(d,\gamma,\theta,\delta_0,K)>0$
 so that  the assertion of the theorem holds if
$\kappa\leq \kappa_0$. Fix $x_0\in \cO$ and denote
$a^{ij}_0(t,x)=a^{ij}(t,x_0)$ and $\sigma^{ik}_0(t,x)=\sigma^{ik}(t,x_0)$. Then $u$ satisfies
$$
du=(a^{ij}_0u_{x^ix^j}+
 f_0)\,dt
+(\sigma^{ik}_0 u_{x^i}+g^{k}_0)\,dw^{k}_{t}, \quad u(0)=u_0,
$$
where
$$
f_0=(a^{ij}-a^{ij}_0) u_{x^ix^j}+b^iu_{x^i}+ cu+f, \quad g^{ik}_0=(\sigma^{ik}-\sigma^{ik}_0)u_{x^i}+\mu^k u+g^k.
$$
 By Corollary \ref{corollary good},
\begin{equation}
                            \label{good afternoon}
\|u\|_{\frH^{\gamma+2}_{p,\theta}(\cO,T)}\leq
N\left(\|\psi f_0\|_{\bH^{\gamma}_{p,\theta}(\cO,T)}
+\|g_0\|_{\bH^{\gamma+1}_{p,\theta}(\cO,T,\ell_2)}
+\|u_0\|_{U^{\gamma+2}_{p,\theta}(\cO)}\right).
\end{equation}

 If $\gamma$ is not integer,  then by Lemma \ref{lemma 2.14}(iii) with  some $\nu \in (0,1-\frac{\gamma}{\gamma_+})$ (e.g. $\nu=\frac{1}{2}(1-\frac{\gamma}{\gamma_+})$),
$$
\|(a^{ij}-a^{ij}_0)\psi
u_{x^ix^j}\|_{\bH^{\gamma}_{p,\theta}(\cO,T)}\leq N
\sup|a^{ij}-a^{ij}_0|^{\nu}\|\psi
u_{x^ix^j}\|_{\bH^{\gamma}_{p,\theta}(\Omega,T)}
\leq
N\kappa^{\nu}\|\psi^{-1}u\|_{\bH^{\gamma+2}_{p,\theta}(\Omega,T)},
$$
$$
\|\psi b^iu_{x^i}+\psi cu\|_{\bH^{\gamma}_{p,\theta}(\cO,T)}\leq N\sup |\psi b^i|^{\nu} \|u_x\|_{\bH^{\gamma}_{p,\theta}(\cO,T)}+N\sup |\psi^2 c|\|\psi^{-1}u\|_{\bH^{\gamma}_{p,\theta}(\cO,T)},
\leq N\kappa^{\nu}\|\psi^{-1}u\|_{\bH^{\gamma+2}_{p,\theta}(\Omega,T)},
$$
 and similarly
$$
\|(\sigma^{i}-\sigma^i_0)u_x\|_{\bH^{\gamma+1}_{p,\theta}(\cO,T,\ell_2)}\leq N \sup |\sigma^i-\sigma^i_0|^{\nu}_{\ell_2}\|u_x\|_{\bH^{\gamma+1}_{p,\theta}(\cO,T)}\leq N\kappa^{\nu}\|\psi^{-1}u\|_{\bH^{\gamma+2}_{p,\theta}(\cO,T)},
$$
$$
\|\mu^k u\|_{\bH^{\gamma+1}_{p,\theta}(\cO,T,\ell_2)}\leq N\sup |\psi \mu|^{\nu}_{\ell_2} \|\psi^{-1}u\|_{\bH^{\gamma+1}_{p,\theta}(\cO,T)}\leq N\kappa^{\nu}\|\psi^{-1}u\|_{\bH^{\gamma+2}_{p,\theta}(\cO,T)}.
$$
By these and (\ref{good afternoon}),
\begin{equation}
                          \label{eqn 8.8.8}
\|u\|_{\frH^{\gamma+2}_{p,\theta}(\cO,T)}\leq  N\kappa^{\nu}\|\psi^{-1}u\|_{\bH^{\gamma+2}_{p,\theta}(\cO,T)}+
N\left(\|\psi f\|_{\bH^{\gamma}_{p,\theta}(\cO,T)}
+\|g\|_{\bH^{\gamma+1}_{p,\theta}(\cO,T,\ell_2)}
+\|u_0\|_{U^{\gamma+2}_{p,\theta}(\cO)}\right).
\end{equation}
Thus it is enough to take $\kappa_0$ so that $N\kappa^{\nu}<1/2$
for all $\kappa\leq \kappa_0$.

If $\gamma=0$, then obviously
\begin{eqnarray*}
&&\|\psi (a^{ij}-a^{ij}_0)u_{x^ix^j} +\psi b^iu_{x^i}+\psi cu\|_{\bL_{p,\theta}(\cO,T)}\\&\leq& \sup|a^{ij}-a^{ij}_0|\|\psi u_{xx}\|_{\bL_{p,\theta}(\cO,T)}
+\sup |\psi b|\|u_x\|_{\bL_{p,\theta}(\cO,T)}+\sup|\psi^2 c| \|\psi^{-1}u\|_{\bL_{p,\theta}(\cO,T)}\\
&\leq& N\kappa \|\psi^{-1}u\|_{\bH^2_{p,\theta}(\cO,T)},
\end{eqnarray*}
and by Lemma \ref{lemma 2.14} (also see (\ref{see this})) with $\nu=\varepsilon/(1+\varepsilon)$,
$$
\|(\sigma^i-\sigma^i_0)u_x\|_{\bH^1_{p,\theta}(\cO,T)}\leq N|\sigma^i-\sigma^i_0|^{(0)}_1 \|u_x\|_{\bH^1_{p,\theta}(\cO,T)}\leq
 N\sup |\sigma^i-\sigma^i_0|^{\nu}\|u_x\|_{\bH^1_{p,\theta}(\cO,T)}\leq N\kappa^{\nu}\|\psi^{-1}u\|_{\bH^2_{p,\theta}(\cO,T)},
 $$
 $$
 \|\mu u\|_{\bH^1_{p,\theta}(\cO,T)}\leq N|\psi \mu|^{(0)}_1 \|\psi^{-1}\|_{\bH^1_{p,\theta}(\cO,T)}\leq
 N\sup |\psi \mu|^{\nu}\|\psi^{-1}u\|_{\bH^1_{p,\theta}(\cO,T)}\leq N\kappa^{\nu}\|\psi^{-1}u\|_{\bH^2_{p,\theta}(\cO,T)}.
$$
These lead to  (\ref{eqn 8.8.8}) for $\gamma=0$.

If $\gamma=1,2,3,...$, then by Lemma \ref{lemma 2.14}(ii)
\begin{eqnarray*}
\|(a^{ij}-a^{ij}_0)\psi
u_{x^ix^j}\|_{\bH^{\gamma}_{p,\theta}(\Omega,T)}&\leq&
N \sup|a^{ij}-a^{ij}_0|\|\psi
u_{x^ix^j}\|_{\bH^{\gamma}_{p,\theta}(\Omega,T)}+N|a^{ij}|^{(0)}_{\gamma}\|\psi u_{xx}\|_{\bH^{\gamma-1}_{p,\theta}(\cO,T)}\\
&\leq& N\kappa \|\psi^{-1}u\|_{\bH^{\gamma+2}_{p,\theta}(\cO,T)}+N\|\psi^{-1}u\|_{\bH^{\gamma+1}_{p,\theta}(\cO,T)},
\end{eqnarray*}
and similarly,
$$
\|\psi b^iu_{x^i}+\psi cu\|_{\bH^{\gamma}_{p,\theta}(\cO,T)} \leq N\kappa \|\psi^{-1}u\|_{\bH^{\gamma+2}_{p,\theta}(\cO,T)}+N\|\psi^{-1}u\|_{\bH^{\gamma+1}_{p,\theta}(\cO,T)},
$$
$$
|(\sigma^i-\sigma^i_0)u_x\|_{\bH^{\gamma+1}_{p,\theta}(\cO,T,\ell_2)}\leq N\kappa \|\psi^{-1}u\|_{\bH^{\gamma+2}_{p,\theta}(\cO,T)}+N\|\psi^{-1}u\|_{\bH^{\gamma+1}_{p,\theta}(\cO,T)}.
$$
Thus if $\kappa_1$ is sufficiently small and $\kappa\leq \kappa_1$, then
\begin{eqnarray*}
\|u\|_{\frH^{\gamma+2}_{p,\theta}(\cO,T)}\leq N \|\psi^{-1}u\|_{\bH^{\gamma+1}_{p,\theta}(\cO,T)}
+ N\left(\|\psi f\|_{\bH^{\gamma}_{p,\theta}(\cO,T)}
+\|g\|_{\bH^{\gamma+1}_{p,\theta}(\cO,T,\ell_2)}
+\|u_0\|_{U^{\gamma+2}_{p,\theta}(\cO)}\right).
\end{eqnarray*}
This and  the inequality
$$
\|\psi^{-1}u\|_{H^{\gamma+1}_{p,\theta}(\cO)}\leq \varepsilon \|\psi^{-1}u\|_{H^{\gamma+2}_{p,\theta}(\cO)}+N(\varepsilon)\|\psi^{-1}u\|_{H^{2}_{p,\theta}(\cO)},
$$
yield
$$
\|u\|_{\frH^{\gamma+2}_{p,\theta}(\cO,T)}\leq N \|\psi^{-1}u\|_{\bH^{2}_{p,\theta}(\cO,T)}
+ N\left(\|\psi f\|_{\bH^{\gamma}_{p,\theta}(\cO,T)}
+\|g\|_{\bH^{\gamma+1}_{p,\theta}(\cO,T,\ell_2)}
+\|u_0\|_{U^{\gamma+2}_{p,\theta}(\cO)}\right).
$$
Take $\kappa_0=\kappa_0(0)$ chosen in the above when $\gamma=0$.
Then it suffices to take $\kappa_0=\kappa_0(\gamma)$ so that
$\kappa_0< \kappa_0(0)\wedge \kappa_1$.

{\bf{Step 2}}. We generalize the result of Step 1 by summing up the local estimations of $u$.

Let $x_0\in \partial \cO$.  Fix a nonnegative
function $\eta\in C^{\infty}_0(B_1(0))$ so that $\eta(x)=1$ for
$|x|\leq 1/2$ and define $\eta_n(x)=\eta(n(x-x_0))$,
$$
a^{ij}_n(t,x)=a^{ij}(t,x)\eta_n(x)+(1-\eta_n(x))a^{ij}(t,x_0)=a^{ij}(t,x_0)+\eta_n(x)(a^{ij}(t,x)-a^{ij}(t,x_0)),
$$
$$
\sigma^{ik}_n(t,x)=\eta_n(x)\sigma^{ik}(t,x)+(1-\eta_n(x))\sigma^{ik}(t,x_0),
$$
$$
b^i_n=b^i\eta_n,  \quad c_n=c\eta_n, \quad \mu_n=\eta_n \mu.
$$
 Then
 $$
 |a^{ij}_n(t,x)-a^{ij}_n(t,y)|\leq 2 \sup_{x\in supp\, \, \eta_n}\eta_n(x)|a^{ij}(t,x)-a^{ij}(t,x_0)|,
 $$
 $$
 |\sigma^{ij}_n(t,x)-\sigma^{ik}_n(t,y)|_{\ell_2}\leq 2\sup_{x\in supp\, \, \eta_n} \eta_n(x)|\sigma^{ik}(t,x)-\sigma^{ik}(t,x_0)|_{\ell_2}
 $$
 and for any multi-index $\alpha$,
 $$
 \sup_n \sup_{x\in \cO} \psi^{|\alpha|} |D^{\alpha}\eta_n|<N(|\alpha|,\eta)<\infty.
 $$
 Indeed, for instance, if $x$ is in the support of $\eta_n$, then $\rho(x)\leq 1/n$ and thus $|\rho(x)D\eta_n(x)|=n\rho(x)|\eta_x (n(x-x_0))|\leq \sup_x|\eta_x|$.
Using this one can easily   check that the coefficients
$a^{ij}_n, b^i_n, \cdots, \mu^k_n$ satisfy (\ref{1.23.01}), (\ref{eqn 2.20.1}) and (\ref{see this}) with some constant $K_0$,
which is independent of $n$.

 Take $\kappa_0$ from Step 1
corresponding to $d,\gamma,\delta_0,K,K_0$ and $\theta$. We fix $n$
large enough so that
$$
|a^{ij}_n(t,x)-a^{ij}_n(t,y)|+|\sigma^i_n(t,x)-\sigma^i_n(t,y)|_{\ell_2}+|\psi b^i_n(t,x)|+|\psi^2c_n(t,x)|+|\psi \mu_n|_{\ell_2}< \kappa_0
\quad  \forall \omega,t,x,y.
$$
This is possible due to the uniform continuity of $a^{ij}, \sigma^i$ and condition (\ref{eqn 2.20.2}).

Now we denote $v=u\eta_{2n}$. Then since $\eta_n=1$ and e.g. $a^{ij}_n=a^{ij}$ on the support of $v$,  $v$ satisfies
$$
dv=(a^{ij}_nv_{x^ix^j}+b^i_nv_{x^i}+
c_nv+\bar{f})dt+(\sigma^{ik}_nv_{x^i}+\mu^k_nv+\bar{g}^k_n)dw^k_t, \quad v(0)=u_0\eta_{2n},
$$
where
$$
\bar{f}:=-2a^{ij}u_{x^i}\eta_{2nx^j}-a^{ij}u\eta_{2nx^ix^j}-b^iu\eta_{nx^i}+\eta_{2n}f, \quad \bar{g}^k=-\sigma^{ik}u\eta_{2nx^i}+\eta_{2n}g^k.
$$
By the result of Step 1, for each $t\leq T$,
\begin{eqnarray*}
\|v\|^p_{\frH^{\gamma+2}_{p,\theta}(\Omega,t)}&\leq& N(
\|\psi\bar{f}\|^p_{\bH^{\gamma}_{p,\theta}(\cO,t)}
+\|\bar{g}\|^p_{\bH^{\gamma+1}_{p,\theta}(\cO,t)}
+\|u_0\eta_{2n}\|^p_{U^{\gamma+2}_{p,\theta}(\cO)})\\
&\leq &N\|\psi u_x\|^p_{\bH^{\gamma}_{p,\theta}(\cO,t)}+N\|u\|^p_{\bH^{\gamma+1}_{p,\theta}(\cO,t)}+N(\|\psi f\|^p_{\bH^{\gamma}_{p,\theta}(\cO,t)}
+\|g\|^p_{\bH^{\gamma+1}_{p,\theta}(\cO,t,\ell_2)}
+\|u_0\|^p_{U^{\gamma+2}_{p,\theta}(\cO)})\\
&\leq &N\|u\|^p_{\bH^{\gamma+1}_{p,\theta}(\cO,t)}+N(\|\psi f\|^p_{\bH^{\gamma}_{p,\theta}(\cO,t)}
+\|g\|^p_{\bH^{\gamma+1}_{p,\theta}(\cO,t,\ell_2)}
+\|u_0\|^p_{U^{\gamma+2}_{p,\theta}(\cO)}),
\end{eqnarray*}
where $N$ is independent of $t$, and the second inequality is due  to Lemma \ref{lemma 2.14} and the following:
$$
|a^{ij}\eta_{2nx}|^{(0)}_{\gamma_+}+|\psi a^{ij}\eta_{2nxx}|^{(0)}_{\gamma_+}+|\psi b^i \eta_{2nx}|^{(0)}_{\gamma_+}+|\sigma^i\eta_{2nx}|^{(0)}_{(\gamma+1)_+}
\leq N<\infty.
$$

 Now to estimate $u$, one introduces a partition of
unity $\zeta_i, i=0,1,...,N$ (remember we assume $\cO$ is bounded) so that $\zeta_0\in
C^{\infty}_0(\cO)$ and $\zeta_i=\eta(2n(x-x_i))$, $x_i\in
\partial \cO$ for $i\geq 1$. Then by the above result, for each $i\geq 1$ and $t\leq T$,
\begin{equation}
                  \label{eqn 11111}
\|\zeta_i u\|^p_{\frH^{\gamma+2}_{p,\theta}(\Omega,t)}
\leq N(\|u\|^p_{\bH^{\gamma+1}_{p,\theta}(\cO,t)}+\|\psi f\|^p_{\bH^{\gamma}_{p,\theta}(\cO,t)}
+\|g\|^p_{\bH^{\gamma+1}_{p,\theta}(\cO,t,\ell_2)}
+\|u_0\|^p_{U^{\gamma+2}_{p,\theta}(\cO)}).
\end{equation}
Note that since $\zeta_0$ has compact support in $\cO$, for any $h\in H^{\gamma}_{p,\theta}(\cO)$ we have $\zeta_0 h\in H^{\gamma}_p$. Moreover for any $\nu\in \bR$,
\begin{equation}
                           \label{eqn 8.8.9}
\|\psi^{\nu}\zeta_0 h\|_{H^{\gamma}_{p,\theta}(\cO)} \sim \|\psi^{\nu} \zeta_0 h\|_{H^{\gamma}_p}\sim \|\zeta_0 h\|_{H^{\gamma}_p}.
\end{equation}
Write  down the equation for $\zeta_0 u$ and apply  Theorem 5.1 of \cite{Kr99} to get
\begin{eqnarray*}
\|\zeta_0u\|^p_{\frH^{\gamma+2}_{p,\theta}(\cO,t)} \sim \|\zeta_0u\|^p_{\cH^{\gamma+2}_{p}(t)}
&\leq&N\|-2a^{ij}u_x\zeta_{0x}-a^{ij}u\zeta_{0xx}-b^iu\zeta_{0x}+\zeta_0f\|^p_{\bH^{\gamma}_p(t)}\\
&&+N\|\sigma^iu\zeta_{0x^i}+\zeta_0g\|^p_{\bH^{\gamma+1}_p(t)}+N\|\zeta_0u_0\|^p_{U^{\gamma+2}_p}.
\end{eqnarray*}
Actually the smoothness condition on the coefficients in Theorem 5.1 of \cite{Kr99} is different from ours since there the coefficients are assumed to be in standard
H\"older spaces. But since $\zeta_0$ has compact support, one can replace these coefficients with $\bar{a}^{ij},\bar{b}^i,\cdots, \bar{\mu}^k$ having finite standard H\"older norms without hurting the equation. By (\ref{eqn 8.8.9}),
$$
\|\bar{a}^{ij}u_x\zeta_{0x}\|_{\bH^{\gamma}_p(t)} \leq N\|u_x\zeta_{0x}\|_{\bH^{\gamma}_p(t)}\leq
N\|\psi u_x \zeta_{0x}\|_{\bH^{\gamma}_{p,\theta}(\cO,t)}\leq N \|\psi u_x \|_{\bH^{\gamma}_{p,\theta}(\cO,t)}\leq N\|u\|_{\bH^{\gamma+1}_{p,\theta}(\cO,t)}.
$$
Similar calculus easily shows $\zeta_0 u$ also satisfies (\ref{eqn 11111}). By summing all these estimates and using (\ref{gronwall}) we get, for $t\leq T$
\begin{eqnarray*}
\|u\|^p_{\frH^{\gamma+2}_{p,\theta}(\cO,t)}&\leq& N \|u\|^p_{\bH^{\gamma+1}_{p,\theta}(\cO,t)}
+N\|\psi f\|^p_{\bH^{\gamma}_{p,\theta}(\cO,t)}+N\|g\|^p_{\bH^{\gamma+1}_{p,\theta}(\cO,t)}
+N\|u_0\|^p_{U^{\gamma+2}_{p,\theta}(\cO)}\\
&\leq& N \int^t_0 \|u\|^p_{\frH^{\gamma+2}_{p,\theta}(\cO,s)} ds+ N \left(\|\psi f\|^p_{\bH^{\gamma}_{p,\theta}(\cO,T)}
+\|g\|^p_{\bH^{\gamma+1}_{p,\theta}(\cO,T)}+
\|u_0\|_{U^{\gamma+2}_{p,\theta}(\cO)}\right).
\end{eqnarray*}
 Thus estimate (\ref{eqn 08.9.1})  follows from this
and Gronwall's inequality.
\qed

\mysection{Proof of Theorem \ref{theorem 1}}
                                                 \label{section theorem 1}
Due to  the method of continuity and a priori estimate (\ref{eqn 08.9.1}) (see e.g. the
proof of Theorem 5.1 of \cite{Kr99} for details), to finish the proof, we only show that for any $f\in
\psi^{-1}\bH^{\gamma}_{p,\theta}(\cO,T),g\in\bH^{\gamma+1}_{p,\theta}(\cO,T)$ and $u_0\in U^{\gamma+2}_{p,\theta}(\cO)$,
the equation
\begin{equation}
                           \label{3.04.02}
du=(\Delta u+f)\,dt+ g^k dw^k_t, \quad u(0)=0
\end{equation}
has a solution $u\in \frH^{\gamma+2}_{p,\theta}(\cO,T)$. We can
approximate $g=(g^1,g^2,...)$ with functions
 having only finite nonzero entries, and
  smooth functions with compact support are dense in $H^{\nu}_{p,\theta}(\cO)$.
  Therefore it follows from a priori estimate (\ref{eqn 08.9.1}) that, to prove existence of solution, we
 may assume that $g$ has only finite nonzero entries and
 is bounded on $\Omega \times [0,T] \times \cO$ along
with each derivative in $x$ and vanishes if $x$ is near $\partial
\cO$. Indeed, let $g^n \to g$ in $\bH^{\gamma+1}_{p,\theta}(\cO,T,\ell_2)$ where $g^n$ satisfy the above mentioned conditions, and assume that
equation (\ref{3.04.02}) with $g^n$ in place of $g$ has a solution $u^n$, then using  (\ref{eqn 08.9.1}) applied for $u^n-u^m$ one easily finds that
$\{u^n\}$ is a Cauchy sequence in $\frH^{\gamma+2}_{p,\theta}(\cO,T)$ and $u_n \to u$ for some $u\in \frH^{\gamma+2}_{p,\theta}(\cO,T)$. Obviously the limit $u$ becomes a solution of (\ref{3.04.02}) (see Theorem \ref{thm cauchy}).

Under such assumed conditions on $g$,
$$
v(t,x):=\int^t_0 g^k(s,x) dw^k_s
$$
is infinitely differentiable in $x$ and vanishes near $\partial
\cO$. Therefore we conclude $v\in \frH^{\nu}_{p,\theta}(\cO,T)$ for
any $\nu \in \bR$.
 Observe that equation (\ref{3.04.02}) can be written as
$$
d\bar{u}=(\Delta \bar{u} +f+ \Delta v)dt,
$$
where $\bar{u}:=u-v$. Thus we reduced the case to the case in which
$g\equiv 0$.   The same argument shows that we may further assume that $f, u_0$
are bounded along each derivative in $(t,x)$ and vanish near
$\partial \cO$.  Furthermore by   considering $u-u_0$, we find that we also may assume $u_0=0$.

First, we consider the case $\theta\geq d-2+p$.
\begin{lemma}
                             \label{lemma 2.12.1}
Let $\theta \geq d-2+p$, $f\in \bL_{p,d}(\cO,T)$ vanish near $\partial \cO$,
say $f(t,x)=0$ for $x\not\in \cO_k:=\{x\in \cO: \psi(x)>1/k\}$ for some $k>0$. Also assume that
the first derivatives of $f$ in $x$ exist and are bounded. Then the
equation
\begin{equation}
                        \label{eqn 2.12}
du=(\Delta u+f)\,dt, \quad u(0)=0
\end{equation}
 has a solution $u\in
\frH^1_{p,\theta}(\cO,T)$.
\end{lemma}
\pf
By Lemma \ref{lemma 4.2} we only need to prove that there exists a
solution $u\in \psi\bL_{p,d-2+p}(\cO,T)$. Let $n>k$. Since
$\partial \cO_n \in C^{\infty}$, by Theorem 2.10 in
\cite{KK2} (c.f. Theorem IV 5.2 in \cite{O}), there is a unique
(classical) solution $u^n\in \frH^2_{p,d}(\cO_n,T)$ of
$$
du^n=(\Delta u^n +f)dt, \quad u^n(0,\cdot)=0,
$$
such that  $u^n|_{\partial \cO_n}=0$ and $Du^n, D^2u^n$ are
bounded in $[0,T]\times \cO_n$. Extend $u^n(x)=0$ for $x\not\in
\cO_n$, then $u^n$ is Lipschitz continuous in $\cO$. Since  for any $q\geq 2$, $(|u|^q)_t=q|u|^{q-2}uu_t=q|u|^{q-2}u (\Delta u+f)$, for each $x\in \cO_n$,
$$
|u^n(T,x)|^q=q\int^T_0|u^n|^{q-2}u^n(\Delta u^n+f)dt.
$$
Integrate this over $\cO_n$ and do integration by parts to get
$$
\int^T_0\int_{\cO_n}|u^n|^{q-2}|Du^n|^2\,dx\,dt\leq 1/(q-1)\int^T_0\int_{\cO_n}|\psi^{-1}u^n|^{q-1}|\psi f| \psi^{q-2}dx
$$
\begin{equation}
\leq \varepsilon \int^T_0\int_{\cO}|\psi^{-1}u^n|^q \psi^{q-2}\,dx\,dt
+N(\varepsilon,q)\int^T_0\int_{\cO}|\psi f|^q \psi^{q-2}\,dx\,dt. \label{eqn 8.4.10}
\end{equation}
Taking $q=2$ and using Hardy inequality, we get
$$
\sup_n
(\|\psi^{-1}u^n\|_{\bL_{2,d}(\cO,T)}+\|Du^n\|_{\bL_{2,d}(\cO,T)})<\infty.
$$
Now we choose $\zeta^n\in C^{\infty}_0(\cO_n)$ such that
$\zeta^n=1$ on $\cO_k$,  $ \psi\zeta^n_{x}, \psi^2\zeta^n_{xx}$
are bounded in $\cO$ uniformly in $n$, and $\zeta^n(x) \to 1$ for $x\in \cO$ as
$n\to \infty$. Then $u^n\zeta^n\in \frH^2_{2,d}(\cO,T)$
satisfies
$$
(u^n\zeta^n)_t=\Delta(u^n\zeta^n)-2u^n_{x^i}\zeta^n_{x^i}
-u^n\Delta \zeta^n+f.
$$
By  a priori estimate (\ref{eqn 08.9.1})
$$
\|u^n\zeta^n\|_{\frH^2_{2,d}(\cO,T)}\leq
N\|u^n_{x^i}\psi\zeta^n_{x^i}-\psi^{-1}u^n\psi^2\Delta\zeta^n\|_{\bL_{2,d}(\cO,T)}+N\|\psi f\|_{\bL_{2,d}(\cO,T)}.
$$
 By dominated convergence theorem,
$$
\|u^n_{x^i}\psi\zeta^n_{x^i}-\psi^{-1}u^n\psi^2\Delta\zeta^n\|_{\bL_{2,d}(\cO,T)}
\to 0 \quad \text{as}\quad n \to \infty.
$$
 Denote $v^n=u^n\zeta^n\in
\frH^1_{2,d}(\cO,T)$, then $\{v^n\}$ is a bounded sequence
in $\frH^1_{2,d}(\cO,T)$. By Theorem \ref{thm cauchy} there exists $u\in \frH^1_{2,d}(T)$ so that $v^n$ and  $\bD u^n$ converges weakly to $u$ and $\bD u$ respectively, and for any $\phi\in C^{\infty}_0(\cO)$ and $t\in [0,T]$ we have $(v^n(t),\phi) \to (u(t),\phi)$ weakly in $L_2(\Omega)$. Since $v^n \to u$ weakly in $\bH^1_{2,d-2}(\cO,T)$, we have $\Delta v^n \to v$ in $\bH^{-1}_{2,d+2}(\cO,T)$.  These and  the fact that $(-2u^n_{x^i}\zeta^n_{x^i}
-u^n\zeta^n_{x^ix^j}, \phi)=0$ for all large $n$ show that $u$ satisfies (\ref{eqn 2.12}) in the sense of distribution.

Also,  (\ref{eqn 8.4.10}) with $q=p$  and (\ref{eqn 08.2.1}) certainly show that $\sup_n \|\psi^{-1}u^n\|_{\bL_{p,d-2+p}(\cO,T)}<\infty$. It follows that
$\psi^{-1}u\in \bL_{p,d-2+p}(\cO,T)\subset \bL_{p,\theta}(\cO,T)$. The lemma is proved.
\qed

To finish the proof, we only need to show that there exists $\beta_1>0$ so that $\theta > d-2+p-\beta_1$, then equation (\ref{eqn 2.12}) has a solution $u\in \bL_{p,\theta-p}(\cO,T)$. As before we assume $f$ is sufficiently smooth and vanishes near the boundary. Take $\kappa_0$ from Step 1 of the proof of Theorem \ref{a priori}. We already proved that if $|\psi b^i|+|\psi^2 c|\leq \kappa_0$
and $\theta=d-2+p$, the equation
\begin{equation}
                  \label{eqn last}
dv=(\Delta v+b^iv_{x^i}+cv +\psi^{\beta}f), \quad v(0)=0
\end{equation}
has a unique solution $v\in \frH^1_{p,\theta}(\cO,T)$ for any $\beta$. Since $\psi_x$ and $\psi \psi_{xx}$ are bounded we can fix $\beta>0$ so that for
$$
b^i:=2\psi^{\beta}(\psi^{-\beta})_{x^i}=-2\beta\psi^{-1}\psi_{x^i},
$$
$$
 c:=\psi^{\beta}\Delta (\psi^{-\beta})=\beta(\beta-1)\psi^{-2}|\psi_x|^2-\beta\psi^{-1}\Delta \psi
$$
the inequality $|\psi b^i|+|\psi^2 c|\leq \kappa_0$ holds, and thus (\ref{eqn last}) has a solution $v\in \frH^1_{p,d-2+p}$.
Now it is enough to check that $u:=\psi^{-\beta}v$ satisfies (\ref{eqn 2.12}) and $u\in \frH^1_{p,d-2+p-\beta p}(\cO,T)\subset \frH^1_{p,\theta}(\cO,T)$
for any $\theta\geq d-2+p-\beta p$. The theorem is proved.

\mysection{Proof of Theorem \ref{theorem 22}}
                           \label{section theorem 22}
 Our previous proofs (see e.g. Lemma \ref{lemma 4.25}) show that we only need to consider case $\theta=d$ with
    equation (\ref{eqn 08.8}) having coefficients independent of $x$. First observe that inclusion
    $\frH^{\gamma+2}_{p,d}(\cO,T)\subset \frH^{\gamma+2}_{p,d-2+p}(\cO,T)$ gives  the uniqueness result for  free.
     Also Remark \ref{remark 08.5} shows that we only need to show there is a   solution $u\in \bL_{p,d-p}(\cO,T)$, so that
    $$
    \|\psi^{-1}u\|_{\bL_{p,d}(\cO,T)}\leq N\left(\|\psi f\|_{\bL_{p,d}(\cO,T)}+\|g\|_{\bH^1_{p,d}(\cO,T,\ell_2)}+\|u_0\|_{U^2_{p,d}(\cO)}\right).
    $$

    For simplicity, assume $u_0=0$. Denote
    $$
    \cF_{p,\theta}=\{(f,g): \|(f,g)\|_{\cF_{p,\theta}}=\|\psi f\|_{\bL_{p,\theta}(\cO,T)}+\|g\|_{\bH^1_{p,\theta}(\cO,T,\ell_2)}<\infty \}.
    $$
    Fix $q>2$  and $\beta\in (0,\beta_0)$, where $\beta_0=\beta_0(d,\delta_0,K)$. Then by Theorem \ref{theorem 1}, the map
    $\cR: (f,g)\to \psi^{-1}u$, where $u$ is the solution of  equation (\ref{eqn 08.8}) is a bounded operator from $\cF_{2,d-\beta}$ to $\bL_{2,d-\beta}(\cO,T)$, and
    from $\cF_{q,d-2+q}$ to $\bL_{q,d-2+q}(\cO,T)$.  Choose $\nu\in (0,1)$ and  $p\in (2,q)$  so that $d=(1-\nu)(d-\beta)+\nu (d-2+q)$ and  $1/p=(1-\nu)/2+\nu/q$.
    Then $F_{p,d}$ (resp. $\bL_{p,d}(\cO,T)$) becomes a complex interpolation space of $F_{2,d-\beta}$ and $F_{q,d-2+q}$ (resp. $\bL_{2,d-\beta}(\cO,T)$ and $\bL_{q,d-2+q}(\cO,T)$), that is,
    $$
    \cF_{p,d}=[F_{2,d-\beta},F_{q,d-2+q}]_{\nu}, \quad \bL_{p,d}(\cO,T)=[\bL_{2,d-\beta}(\cO,T),\bL_{q,d-2+q}(\cO,T)]_{\nu}.
    $$
   (See Proposition 2.4 of \cite{Lo2} and
Theorem 5.1.2 of \cite{BL} for details).
 It follows from the interpolation theory that $\cR$ is a bounded linear map from $\cF_{p,d}$ to $\bL_{p,d}(\cO,T)$ (see Theorem (a) on Page 59 of \cite{T}).  This proves the claim for above fixed $p$. Now for $2\leq p'\leq p$,
    it is enough to notice that for $\nu'$ so that $1/p^{'}=(1-\nu')/2+\nu'/p$,
    $$
    \cF_{p',d}=[F_{2,d},F_{p,d}]_{\nu'}, \quad \bL_{p',d}(\cO,T)=[\bL_{2,d}(\cO,T),\bL_{p,d}(\cO,T)]_{\nu'}.
    $$
    It follows that $\cR$ is a bounded linear map from $\cF_{p',d}$ to $\bL_{p',d}(\cO,T)$.
The theorem is proved.

\end{document}